\documentclass[reqno,11pt]{amsart}
\usepackage{amsmath, latexsym, amsfonts, amssymb, amsthm, amscd}
\numberwithin{equation}{section}

\newtheorem{theorem}{Theorem}[section]
\newtheorem{lemma}[theorem]{Lemma}
\newtheorem{proposition}[theorem]{Proposition}
\newtheorem{rem}[theorem]{Remark}
\newtheorem{definition}[theorem]{Definition}

\newcommand{\cD}{{\ensuremath{\mathcal D}} }

\newcommand{\cL}{{\ensuremath{\mathcal L}} }

\newcommand{\gb}{\beta}
\newcommand{\gd}{\delta}

\newcommand{\gep}{\varepsilon}       % \ge already exists...

\newcommand{\gs}{\sigma}

\renewcommand{\tilde}{\widetilde}          % wider `tilde'
\DeclareMathSymbol{\leqslant}{\mathalpha}{AMSa}{"36} % nicer `smaller or equal'
\DeclareMathSymbol{\geqslant}{\mathalpha}{AMSa}{"3E} % nicer `larger or equal'
\DeclareMathSymbol{\eset}{\mathalpha}{AMSb}{"3F}     % nicer `emptyset'
\newcommand{\N}{\mathbb{N}}

\DeclareMathOperator{\sign}{sign}

\newcommand{\ind}{\bs{1}}

\def\bs{\boldsymbol}

\title{Copolymer at selective interfaces and pinning potentials: weak coupling limits}

\author{Nicolas Petrelis}

\address{Eurandom, P.O. Box 513, 5600 MB Eindhoven, The Netherlands.}

\address{Laboratoire de Math\'ematiques Raphael Salem, site du Madrillet, 76801 Saint-Etienne du Rouvray, France
CNRS U.M.R. 6085.}

\email{petrelis\@@eurandom.tue.nl}

\date{\today}

\begin{document}

\begin{abstract}
We consider a simple random walk of length $N$, denoted by $(S_{i})_{i\in \{1,\dots,N\}}$, and
we define $(w_i)_{i\geq 1}$ a sequence of centered i.i.d. random variables. For $K\in\N$ we
define $((\gamma_i^{-K},\dots,\gamma_i^K))_{i\geq 1}$ an i.i.d sequence of random vectors.
We
set $\beta\in \mathbb{R}$, $\lambda\geq 0$ and $h\geq 0$, and transform the
measure on the set of
random walk trajectories with the Hamiltonian
$\lambda \sum_{i=1}^{N} (w_i+h) \sign(S_i)+\beta \sum_{j=-K}^{K}\sum_{i=1}^{N} \gamma_{i}^{j}\
\boldsymbol{1}_{\{S_{i}=j\}}$.
This transformed path measure describes an hydrophobic(philic) copolymer
interacting with a layer
of width $2K$ around an interface between oil and water.

In the present article we prove the convergence in the limit of weak coupling
(when $\lambda$, $h$ and $\beta$ tend to $0$)
of this discrete model towards its continuous counterpart.
To that aim we further develop a technique of coarse graining introduced by Bolthausen and den Hollander in \cite{BDH}.
Our result shows, in particular, that the randomness of the pinning around the interface vanishes
as the coupling becomes weaker.

\bigskip


\bigskip
\end{abstract}

\maketitle

\thispagestyle{empty}
%\chapter{Copolymer Pinned at an Interface}\label{chap2}
%In this chapter we consider a copolymer, composed by hydrophobic and hydrophilic monomers.

\section{Introduction and main results}

\subsection{A discrete model of copolymer with adsorption}

We consider a copolymer of $N$ monomers, and an interface separating two solvents (for example oil and water).
The interface runs along the $x$-axis.
The possible configurations of the polymer are given by the trajectories of a simple random walk
$S=(S_i)_{i\geq 1}$ of length $N$ such that $S_0=0$ and $(S_i-S_{i-1})_{i\geq 1}$ is an i.i.d. sequence
of Bernoulli trials satisfying $P(S_1=\pm 1)=1/2$. We let $\Lambda_{i}=\sign(S_{i})$ when $S_i\neq 0$ and
$\Lambda_i=\Lambda_{i-1}$ otherwise.
In size $N$ we take into account the interactions between the polymer and the medium by associating
with each trajectory $S$ the Hamiltonian
\begin{equation}\label{Ham1}
H_{N,\beta,\lambda,h}^{w,\gamma}(S)=\lambda \sum_{i=1}^{N} (w_{i}+h) \Lambda_{i}+
\beta \sum_{j=-K}^{K}\sum_{i=1}^{N} \gamma_{i}^{j}\ \ind_{\{S_{i}=j\}},
\end{equation}
where $\lambda,h\geq 0$, $\beta\in\mathbb{R}$,
$w=(w_i)_{i\geq 1}$ is an i.i.d. sequence of bounded and symmetric
random variables and
$\gamma=((\gamma^{-K}_{i},\dots,\gamma^K_i))_{i\geq 1}$
is an i.i.d. sequence of random vectors ($w$ and $\gamma$ being independent). We stress that
$w$ and $\gamma$ are defined under the probability $\mathbb{P}$ and that
the variables $\gamma_1^{-K},\dots,\gamma_1^K$
are independent but can have different laws.
%\begin{itemize}
%%and the associated probability and expectation will be denoted by $P$ and $E$.
%\item  \underline{Interaction polymer-solvents.} This interaction is given by the Hamiltonian
%$H_1(S)=-\lambda \sum_{i=1}^{N} (w_{i}+h) \Lambda_{i}$, with $\lambda \geq 0$, $h\geq 0$,
%$\Lambda_{i}=\sign(S_{i})$ and $w=(w_i)_{i\geq 1}$ a sequence of i.i.d. random variables.\\
%\item  \underline{Interaction polymer-interface.} This interaction, which involves a pinning potential in a layer
%of finite width around
%the interface is given by the Hamiltonian
%$H_2(S)=\beta \sum_{j=-K}^{K}\sum_{i=1}^{N} \gamma_{i}^{j}\ \ind_{\{S_{i}=j\}}$ where $\beta\in\mathbb{R}$
%and $\gamma=\big(\gamma^{j}_{i}\big)_{i\geq 1,\j\in\{-K,\dots,K\}}$
%is a field of independent random variables.
%\end{itemize}
 This Hamiltonian allows to define
the polymer measure $P_{N,\beta,\lambda,h}^{w,\gamma}$ as
\begin{equation}\label{Z}
\frac{dP_{N,\beta,\lambda,h}^{w,\gamma}}{dP}(S)=
\frac{\exp\big(H_{N,\beta,\lambda,h}^{w,\gamma}(S)\big)}{Z_{N,\beta,\lambda,h}^{w,\gamma}}
.\end{equation}

%\subsection{Copolymer with adsorption}
This discrete model has already been investigated in
physics (see \cite{JSW} or \cite{SW}) and mathematics (see \cite{GiacTon2}) in the case $K=0$
and under the name copolymer with adsorption. This model
%comes mainly from the fact that
%This model is the generalization of the copolymer near a selective interface model, introduced in
%\cite{BDH}, in which a random interaction with a layer of finite width around the interface is added.
is very natural, because it interpolates between two classes of models that
have received a lot of attention in the literature:

\begin{itemize}
\item The \underline{pure pinning} model, which is obtained by setting
$\lambda=0$. In this case only the interaction
with the layer around the origin is activated. This model has been studied
in the case $K=0$, for instance in
\cite{SidAlex}, \cite{Alex}, \cite{JRV}, \cite{Petr}.

\item The \underline{random copolymer} model, which is obtained by fixing
$\beta=0$. In this case only the interaction between the monomers and the two
solvents is activated. It has been studied for instance in
\cite{BDH}, \cite{BodGiac}, \cite{BiskDH}.
\end{itemize}
%Apart from the pinning case with very negative variables $\beta \gamma_i^j$,
In general, these
two models undergo a localization-delocalization phase transition, which results from an energy-entropy
competition. In fact, in both cases, some trajectories are energetically favored with respect to the others.
In the pinning case, it concerns the trajectories that remain close to the interface
to touch the sites that carry a positive reward $\beta \gamma_i^j$ as often as possible. In the
copolymer case, for every $i\in\{1,\dots,N\}$ $w_i+h>0$ (respectively $w_i+h<0$) means that
the $i$-th monomer is hydrophobic (resp. hydrophilic) and therefore, the energetically favored
trajectories cross the interface often to put as many monomers as possible
in their preferred solvent.
In both cases, these favored trajectories are localized in the
neighborhood of the interface. Therefore, they carry much less entropy than the trajectories which
wander away far from the interface.
%\begin{rem}
%\rm{In the single interface case, we do not assume that the rewards at the origin, given by the variables
%$\gamma$, are non negative. On the contrary, in the multi-interfaces model, we assume that the variables $q$
%are non negative. This is essentially to raise the clarity of the article, because with techniques
%similar to those we use in the single-interface model, we could handle the case
%$q\in\mathbb{R}$ and obtain the same results without difficulty.}
%\end{rem}

At this stage we introduce the free energy of the system that will be a key tool to define
the localized and delocalized regimes. Thus, for $N\in\mathbb{N}$ and every disorder $(w,\gamma)$
we define
$\Phi^{w,\gamma}_{N}$ as
\begin{equation}
\frac{1}{N} \log Z_{N,\beta,\lambda,h}^{w,\gamma}=\Phi^{w,\gamma}_{N}(\beta,\lambda,h).
\end{equation}
We recall that $(w,\gamma)$ are defined under the law $\mathbb{P}$ and
we denote by $\Phi_{N}(\beta,\lambda,h)$ the quantity
$\mathbb{E}(\Phi^{w,\gamma}_{N}(\beta,\lambda,h))$. Henceforth, we assume that
$\mathbb{E}(\exp(\beta|\gamma^{j}_{1}|))<\infty$ for every $\beta \in\mathbb{R}$ and
$j\in\{-K,\dots,K\}$.

\begin{proposition}\label{prop1}
For every $\beta\in\mathbb{R}$, $\lambda \geq 0$, $h\geq 0$, there exists a non random real number, denoted by
$\Phi(\beta,\lambda,h)$, such that $\mathbb{P}$ almost surely in $(w,\gamma)$ %$\mathbb{P}$ a.s. in $\chi$ the following convergence,
\begin{equation*}
\lim_{N\to \infty} \Phi_{N}^{w,\gamma}(\beta,\lambda,h)=\Phi(\beta,\lambda,h).
\end{equation*}
This convergence occurs also in $\mathbb{L}^1$, which entails the convergence of
$\Phi_{N}(\beta,\lambda,h)$ to $\Phi(\beta,\lambda,h)$ as $N$ tends to $\infty$.
The limit $\Phi(\beta,\lambda,h)$ is called the free energy of the model.
\end{proposition}
\noindent This proposition has been proven in different
papers for quantities similar to
$Z_{N,\beta,\lambda,h}^{w,\gamma}$ (see \cite{Giac} or \cite{Giac2} for example).
In our case, the difference comes from the fact that the disorder is spread out over a layer
of finite width  around the interface, but the proof remains essentially the same and is left to the reader.
We also notice that
$\Phi(\beta,\lambda,h)$ is continuous and separately convex.

\subsection{The continuous model}% \small({proposition $2$})}

We define in this section the continuous counterpart of the discrete model.
%In this section, we consider a polymer of length $t\in \mathbb{R}$. The parameters
%$\lambda,h$ and $\beta$ are still non negative.
In size $t$, the configurations of the polymer are given by the set of trajectories
of the Brownian motion $\left(B_{s}\right)_{s\in [0,t]}$.
The Hamiltonian associated with every trajectory $B$ is
\begin{equation}
\tilde{H}_{\beta,\lambda,h}^{R,t}(B)=\lambda \int_{0}^{t} \Lambda(s) (dR_{s}+h ds)+\beta L^0_{t},
\end{equation}
where $L^0_t$ (or $L_{t}$ when there is no ambiguity) is the local time spent at $0$
by $B$ between time $0$ and time $t$.
As in the discrete case we set $\lambda,h\geq 0$, $\beta\in \mathbb{R}$ and $\Lambda_{s}=\sign(B_{s})$.
We denote by $\tilde{\mathbb{P}}$ the law of $R=\left(R_{s}\right)_{s\geq 0}$, which is a standard Brownian
motion, independent of $B$ such that $dR_s$ plays the role of $w_i$.

%We consider $dR_{s}$ an elementary variation of $R$ at position $s$.
%This quantity gives the hydrophobicity of the polymer around the position $s$, and plays
%the role of $w_{i}$ in the discrete model.

As in the discrete case, we define the polymer measure
of length $t$ by perturbing the
law $\tilde{P}$ of the Brownian motion $B$ as follows
\begin{equation}\label{Z}
\frac{d\tilde{P}_{\beta,\lambda,h}^{R,t}}{d\tilde{P}}(B)=
\frac{\exp\big(\tilde{H}_{\beta,\lambda,h}^{R,t}(B)\big)}{\tilde{Z}_{\beta,\lambda,h}^{R,t}}
.\end{equation}
For every $t>0$ and every disorder $R$ we introduce the free energy of the system of size $t$,
denoted by
$\tilde{\Phi}^{R}_{t}$, as
\begin{equation}
\frac{1}{t} \log \tilde{Z}_{\beta,\lambda,h}^{t,R}=\tilde{\Phi}^{R}_{t}(\beta,\lambda,h).
\end{equation}
We also denote by $\tilde{\Phi}_{t}(\beta,\lambda,h)$ the quantity
$\tilde{\mathbb{E}}(\tilde{\Phi}^{R}_{t}(\beta,\lambda,h))$.

\begin{proposition}\label{prop2}
For every $\beta\in\mathbb{R}$, $\lambda \geq 0$, $h\geq 0$, there exists a non random real number, denoted by
$\tilde{\Phi}(\beta,\lambda,h)$, such that $\tilde{\mathbb{P}}$ almost surely in $R$ %$\mathbb{P}$ a.s. in $\chi$ the following convergence,
\begin{equation*}
\lim_{t\to \infty} \tilde{\Phi}_{t}^{R}(\beta,\lambda,h)=\tilde{\Phi}(\beta,\lambda,h).
\end{equation*}
As in the discrete case this convergence occurs also in $\mathbb{L}^1$, and therefore
$\tilde{\Phi}(\beta,\lambda,h)$,
which is the free energy of the model, is the limit of $\tilde{\Phi}_{t}(\beta,\lambda,h)$ as t tends to
$\infty$.
\end{proposition}
\noindent A proof of Proposition \ref{prop2} in the case $\beta=0$
is available in \cite{Giac}. This proof
is adapted in
\cite{Petr2} to cover the case $\beta\neq 0$.
We also notice that $\tilde{\Phi}(\beta,\lambda,h)$ is continuous, separately convex and non-decreasing
in $\beta$.

\subsection{Localized and delocalized regimes}

In the discrete and the continuous model, the free energy gives us a tool to decide,
for every $(\beta,\lambda,h)$,
whether the system is localized
or not. Observe that if we set $D_{N}=\{ S:S_{i}>K \;\forall \; i \;\in \{K+1,\dots,N\}\}$ and use
$P(D_{N})=(1+o(1)) c/\sqrt{N}$ and the law of large numbers we have $\mathbb{P}$-a.s.
\begin{equation}\label{gloups}
\Phi(\beta,\lambda,h)\geq \liminf_{N \to \infty}\frac{1}{N}\log
\textstyle{E}\Big[\exp\big(\lambda \sum_{i=1}^{N}(w_{i}+h)+\beta\sum_{i=1}^{K}\gamma_{i}^{i}\big)
\ind_{\{D_{N}\}}\Big]\geq \lambda h.
%&\geq \lambda h+\lim_{N \to \infty}\frac{\lambda \sum_{i=1}^{N}w_{i}}{N}
%+\lim_{N \to \infty}\frac{\beta \sum_{i=1}^{K}\gamma_{i}^{i}}{N}+
%\lim_{N \to \infty}\frac{\log\left(P\left(D_{N}\right)\right)}{N}
%\geq \lambda h.
\end{equation}
%The first inferior limit of \eqref{gloups} tends to $0$ when $N$ tends to $\infty$ because the
%law of large number can be applied to $(w_{i})_{i\geq 1}$, the second inferior limit tends to $0$
%because it is a constant ($\sum_{i=1}^{K}\gamma_{i}^{i}$) divided by $N$ and the last one
%tends to $0$ because $P(D_{N})=(1+o(1)) c/\sqrt{N}$ as $N \uparrow\infty$.
We will say that the polymer is delocalized when $\Phi(\beta,\lambda,h)=\lambda h$,
because the trajectories in $D_N$ essentially determine the free energy, and localized
when $\Phi(\beta,\lambda,h)>\lambda h$. The $(\beta,\lambda,h)$-space is divided into a
\underline{localized phase},
denoted by $\mathcal{L}$, and a \underline{delocalized phase}, denoted by $\mathcal{D}$.
It is now well understood (see in particular \cite{Giac2} and \cite{GiacTon2}) that such a free energy dichotomy does correspond
to sharply different path behaviors.

In the continuous case, by considering the subset $\tilde{D}_{t}=\{B:B_{s}>0\ \forall\; s \in [1,t]\}$,
a computation similar to \eqref{gloups} shows that $\tilde{\Phi}(\beta,\lambda,h)\geq\lambda h$.
Therefore we can use the same dichotomy used in the discrete case to characterize $\mathcal{L}$ and
$\mathcal{D}$.

\subsubsection{Critical curve}\label{dud}

For $\gamma$, $w$, $R$, $K$ and $\beta$ fixed, both for the discrete and continuous models
there exists a critical curve $\lambda \mapsto h_{c}^{\beta}(\lambda)$ ($\tilde{h}_{c}^{\beta}(\lambda)$
in the continuous case),
which divides the $(\lambda,h)$-space
into $\mathcal{L}=\{(\lambda,h) : h<h_c^\beta(\lambda)\}$ and
$\mathcal{D}=\{(\lambda,h) : h\geq h_c^\beta(\lambda)\}$. In fact, by differentiating with respect to $h$
we obtain for every $N\geq 1$ and $t>0$ that $\Phi_N(\beta,\lambda,h)-\lambda h$ and $\tilde{\Phi}_t(\beta,\lambda,h)-\lambda h$
are non increasing in $h$. Therefore $\Phi(\beta,\lambda,h)-\lambda h$ and
$\tilde{\Phi}(\beta,\lambda,h)-\lambda h$ are also non increasing in $h$ and we can simply
define $h_{c}^{\beta}(\lambda)=\inf\{h\geq 0 : \Phi(\beta,\lambda,h)-\lambda h=0\}$ and
$\tilde{h}_{c}^{\beta}(\lambda)=\inf\{h\geq 0 : \tilde{\Phi}(\beta,\lambda,h)-\lambda h=0\}$.

The scaling property of the Brownian motion entails the equality
$\tilde{\Phi}(a\beta,a\lambda,ah)=a^2\tilde{\Phi}(\beta,\lambda,h)$ for every $a\geq 0$. From this it follows that
$\tilde{h}_{c}^{\lambda \beta}(\lambda)=\lambda K_{c}^{\beta}$ with
$K_c^\beta=\inf\{h\geq 0 : \tilde{\Phi}(\beta,1,h)-h=0\}$. Notice that the quantity $K_c^\beta$ can be
viewed as a critical curve in the $(\beta,h)$-plane for $\lambda=1$ and is non-decreasing in $\beta$.
Moreover, since $(\beta,h)\mapsto \tilde{\Phi}(\beta,1,h)-h$ is convex, we prove easily that $\beta\mapsto
K_c^\beta$ is convex.

\begin{rem}
\rm{Observe that for some values of $\lambda$ and $\beta$ the critical
value $h_{c}^{\beta}(\lambda)$ can be infinite.
In section \ref{finir}, we give certain conditions under which this happens.
We prove also, for the continuous model, that
$K_c^{\beta}<\infty$ and $\tilde{h}_{c}^{\beta}(\lambda)<\infty$ for all $\beta\in\mathbb{R}$.}
\end{rem}
\noindent As a consequence, $\beta\mapsto K_c^{\beta}$ is continuous on $\mathbb{R}$ because it is
convex and finite.
%As a consequence, for the continuous case, we define $\beta_c=\sup\{\beta : K_c^{\beta}<\infty\}$
%and we can assert that
%$K_{c}^{\beta}$ is finite and continuous on $(-\infty,\beta_c)$ and equal to $\infty$ for $\beta>\beta_c$.
%It seems much harder to investigate the case $\beta=\beta_c$, namely to decide whether $K_{c}^{\beta_c}$
%is finite or not. I will not answer this question in this paper.
%The existence of
%this curve is proven in \cite{BDH}
%for the case $\beta=0$, and can be easily adapted to our case. For
%this reason we will not give the details in this article.
%In \cite{Giac}, the proposition $2$ has been proven in the case $\beta=0$, that is to say without pinning term
%at the interface. In the Appendix A, we give a detailed proof of the case $\beta\geq 0$. To that aim,
%we follow the scheme of the proof exposed in \cite{Giac}, and we modify some steps to take into account
%the presence of the $\beta$ term.
\subsection{Discussion of the model and main results}

Before studying more in depth the mathematical properties of the model, we recall that
one of the physical situations that can be modelled by such systems is a polymer
put in the neighborhood of an interface between two solvents (see \cite{BDH}).
Nevertheless, the models considered up to now do not take into account
that such an interface has a finite width, that is to say, a small layer in which the two solvents are more or less
mixed together. In this sense, the model developed here
gives a more realistic image of an interface. Moreover, this model allows
us to consider other physical situations. For instance, a case
in which micro-emulsions of a third solvent are spread in a thin layer around
the interface.

\subsubsection{\underline{Former results about the model}}

We can roughly classify the results available for polymer models in two categories. On the one hand
the results
concerning the \underline{path behavior} of the polymer. In fact,
the separation between the localized and delocalized phases has an interpretation in terms of trajectories of
the polymer. We refer to \cite{AlbZhou},
\cite{BiskDH}, \cite{Giac2}, \cite{GiacTon2}, \cite{Sinai}
for sharp results about the path behavior of the copolymer in $\mathcal{L}$ and we refer to
\cite{GiacTon} for
further results in $\mathcal{D}$.
On the other hand, the results concerning
the \underline{free energy} ($\Phi$): this problem arises only in $\mathcal{L}$ since $\Phi$ is
constant in $\mathcal{D}$. In this last category, we can mention for instance
the strong results about disordered pinning obtained recently in \cite{Alex}, in particular
concerning the comparison between quenched and annealed critical curve at weak disorder.

 %about the copolymer without adsorption ($\beta=0$) in $\mathcal{D}$.%The result asserts that
%for every family of parameters in $\mathcal{L}$ there exist $c_1,c_2>0$
%such that for every $L\in\mathbb{N}$
%$$\sup_{N\geq 1}\  \sup_{0\leq k\leq N} \mathbb{E}P_{N,w}\left(|S_{k}|\geq L\right)\leq c_1
%\exp(-c_2 L).$$
%For the copolymer with adsorption in the case $K=0$, a sharp estimation is given in \cite{GiacTon2}
%concerning the
%asymptotic behavior of the longest excursion of the polymer. The authors denote by $P_{N,w}$ the
%law of the polymer and by $\Delta_N$ the length of the longest excursion made by the polymer after $N$ steps.
%They prove that for every family of parameters
%in $\mathcal{L}$  there exists a constant $c>0$ such that for every $\gep>0$ the equality
%$$\textstyle \lim_{N\to \infty} P_{N,w}\big(c(1-\gep)<\frac{\Delta_N}{\log N}<c(1+\gep)\big)=1\quad \
%\text{occurs in probability.}$$
%Notice that an expression of the constant $c$ is also given.

%Concerning the path properties in the delocalized phase, we refer to  %The authors prove that in
%$\mathcal{D}$ the time spent by the polymer in the lower half plane after $N$ steps is bounded by
%$O(\log(N))$ and even by $o(1)$ if the parameters are chosen deeply enough in $\mathcal{D}$.

\subsubsection{Regularity and scaling limit of the free energy}
For both the copolymer and the pinning model the free energy $\Phi$ is complicated inside
$\mathcal{L}$ and an important question is to figure out if another phase transition
can occur inside
$\mathcal{L}$.
The answer is partially given for the case $K=0$ in \cite{GiacTon2}, where a proof
of the infinite differentiability of the free energy inside $\mathcal{L}$ is given.
This proof is based on a result, that was first given in \cite{Sinai} and \cite{BiskDH} for
the copolymer without adsorption and
asserts that in $\mathcal{L}$ the laws of the polymer's excursions are exponentially tight.
From this tightness, certain correlation inequalities are deduced that are sufficient to prove the infinite
differentiability of the free energy inside $\mathcal{L}$.
Therefore, there is no other phase transition, at least of finite order, inside the localized phase.

The scaling limit of the discrete model is also a question that has been closely studied recently.
In fact, in the case of the copolymer without adsorption ($\beta=0$),
a continuous model is introduced in \cite{BDH} and it
turns out to be the limit of the discrete model at high temperature, i.e.,
when the coupling parameters $\lambda$ and $h$ tend to $0$. The results in \cite{BDH}
deal with the case of $w$ taking values $\pm 1$ and focus on the free energy, i.e, %and in terms of slope of the critical curve at the origin, i.e.,
\begin{equation}\label{AaA}
\lim_{a\to 0} \frac{1}{a^2}\Phi(0,a\lambda,ah)=\tilde{\Phi}(0,\lambda,h).%\quad\quad\text{and}\quad
%\quad \lim_{\lambda\to 0}\frac{h^0_c(\lambda)}{\lambda}=K_c^0.
\end{equation}
This has been generalized in \cite{GiacTon} to a large class of random variables $w$ and it
is of interest in terms of universality of the Brownian limit as we are going to explain.
Effectively, it shows that, when the
coupling constants become weak, the Brownian models
"attracts" any discrete model, regardless of its charge distribution.
The proof is based on a coarse-graining method. In fact, for fixed parameters $(a\lambda,ah)$
the $N$ steps of the polymer are partitioned into blocks of finite
and constant size $L(a)$. It turns out that the characteristic size of the excursions
for $a$ small is of order $1/a^2$. Then, by choosing $L(a)$ of order $1/a^2$, one can, block by block,
approximate the
free energy per steps of the discrete model by the one of the continuous model.
When $N$ tends to $\infty$ the number of blocks tends to $\infty$,
but an ergodic property of the blocks allows to convert the approximation per block
into the convergence \eqref{AaA} involving the discrete and the continuous free energies
in infinite size. In \cite{BDH} it is shown that
the convergence occurs for the slope of the critical curve at the origin as well, i.e.,
$\lim_{\lambda\to 0}h_c^0(\lambda)/\lambda=K_c^0$.

\subsubsection{\underline{Main results}}

In this article, we extend the scaling limit of the free energy given in \cite{BDH} to the model of a copolymer
with adsorption introduced above. We aim particularly at understanding how the random pinning
is modified at high temperature.
Some zones in the interacting layer around the origin carry a large
number of high rewards and play a particular
role from the localization point of view. Indeed, the chain can target these zones when it goes back
to the origin in order to maximize the rewards.
Consequently, some zones favor the localization of the polymer more than others
(see \cite{SidAlex} and \cite{Petr}).
Here the question is whether the passage to a very weak coupling preserves the randomness of these rewards
or leads to a complete averaging of the disorder.

We answer this question in Theorems \ref{theo3}. In fact, by generalizing the limit
\eqref{AaA} to the case $\beta\neq 0$ we prove the convergence
of the discrete model to
the continuous model, when the parameters tend to $0$ at appropriate speeds.
The associated continuous model has a pinning term at the interface, given by the local time
at $0$ of the Brownian motion $B$. Therefore, the randomness of the pinning term vanishes in
the weak coupling limit.

In what follows, we will use the notation
\begin{equation}
\textstyle{\Sigma=\sum_{j=-K}^{K}\mathbb{E}\big(\gamma_{1}^{j}\big)}.
\end{equation}
With the limit (\ref{corol}) given
in Theorem \ref{theo3}, we prove that
the partial derivatives of $(\lambda,\beta)\mapsto h_c^{\beta}(\lambda)$ at the origin with respect to
any vector $(1,\beta)$ ($\beta\in \mathbb{R}$) are only determined by
the quantity $\beta \Sigma$. This is also an important result in terms
of universality of the continuous limit with respect to the disorder $\gamma$. In fact,
it shows that the shape of
the critical surface close to the origin only depends on $\Sigma$.

Before stating Theorem \ref{theo3}, we recall that the variables $(\gamma_1^j)_{j\in\{-K,\dots,K\}}$
are allowed to have different laws, and we assume
without loss of generality that $\mathbb{E}(w_1^2)=1$.
\begin{theorem}\label{theo3}
Let $\beta\in \mathbb{R}$, $\lambda>0$, $h\geq 0$,  and
$\Sigma=\sum_{j=-K}^{K}\mathbb{E}\big(\gamma_{1}^{j}\big)$. Then

\begin{equation}\label{theo31}
\lim_{a \to 0}\
\frac{1}{a^{2}}\ \Phi(a \beta, a \lambda, a h)=\tilde{\Phi}\left(\beta \Sigma, \lambda ,h\right)
\end{equation}
and
\begin{equation}\label{corol}
\lim_{\delta \to 0} \frac{h_c^{\delta \beta}(\delta)}{\delta}=K_c^{\beta \Sigma}.
\end{equation}
\end{theorem}

We can derive from Theorem \ref{theo3} some relevant information concerning two particular cases of the model.
In a first part we
consider the influence of a deterministic pinning term on the critical curve of a copolymer without adsorption.
In a second part we consider the case of an homopolymer
with adsorption.
\subsection{Two particular cases}
\subsubsection{Influence of a depinning term on the critical curve}
We
consider here
the copolymer model with a deterministic pinning term, i.e., $K=0$, $\gamma_1^0=1$.
Up to now the sensibility of the critical curve $\lambda\mapsto h_c^{0}(\lambda)$ to the presence of
a pinning or depinning term is only very partially understood.
Effectively, in the case $(\lambda,h)\in\mathcal{D}$, one can prove that choosing $\beta$ large enough
is sufficient to obtain $h<h_c^{\beta}(\lambda)$, namely to pass from a delocalized regime to
a localized regime. It can be done for instance by restricting the computation of the
free energy to the random walk trajectories that come back to the origin every
second steps. This leaves open the question whether a small $\beta$ can transform the critical curve.
\vspace{1.5cm}

\setlength{\unitlength}{0.23cm}
\begin{picture}(10,15)(-15,0)

\put(0,0){\vector(1,0){30}}
\put(0,0){\vector(0,1){15}}

\put(-.8,-1.3){$0$}

\put(32,-.3){$\lambda$}

\put(-1,16){$h$}
{\thicklines
   \qbezier(0,0)(5,5)(25,8.5)}
{
   \qbezier(0,0)(6,3)(25,5.5)}

% la legende
\put(-10,17){\small
             Fig.\ 1:}%{\thicklines
\put(9,8){$\mathcal{D}$}
\put(14,4.5){$\mathcal{L}_1$}
\put(19,1.5){$\mathcal{L}_2$}
\put(25.5,10){$h_{c}^{0}(\lambda)$}
\put(26.5,4){$\underline{h}(\lambda)$}
\end{picture}
$\phantom{aaa}$\\

The situation does not get easier when $(\lambda,h)\in\cL$. In this case,
it is useful to divide
$\cL$ into the two regions $\cL_1$ and $\cL_2$ separated by the curve
$\lambda\mapsto \underline{h}(\lambda)=(3/4\lambda) \log \mathbb{E}(\exp(4\lambda w_1/3))$
(see Fig $1$). In fact, the localization
strategy displayed in \cite{BodGiac}
to prove that $\underline{h}(\lambda)\leq h_c^0(\lambda)$ is not sensitive to the presence of a depinning term.
This strategy consists in coming back to the origin only to target rare stretches
of negative $w_i$. These rare stretches are of length $l$, and the energetic contribution of each of them is
of order $l$ whereas
the depinning term contributes an energy
$O(1)$. Thus, for
$h<\underline{h}(\lambda)$ (i.e., $(\lambda,h)\in\cL_2$), we can not chose $\beta<0$ such that
$h\geq h_c^{\beta}(\lambda)$.

The case $(\lambda,h)\in \cL_1$ is harder to investigate and we must recall that the strict inequality
$\underline{h}(\lambda)<h_c^0(\beta)$ is not rigorously proven for the moment. However,
some numerical evidences in \cite{GiacGubCar} shows that $\cL_1$ is not
an empty set and
contrary to what we just said about $\cD$ and $\cL_2$,
the influence of a depinning term in the region $\cL_1$ is not understood at all. This leads
to the following open problem: for $(\lambda,h)\in \mathcal{L}_1$,
namely when $\underline{h}(\lambda)\leq h<h_c^0(\lambda)$,
can we find a large enough depinning term $\beta<0$ that leads to a delocalization, i.e.,
$h\geq h_c^{\beta}(\lambda)$?

From this point of view, Theorem \ref{theo3} is an improvement in the knowledge of
the depinning influence in $\mathcal{L}_1$.
Indeed, even if Theorem \ref{theo3} does not directly answer this open problem, it connects it to another problem
that may be easier to solve. Effectively, if one can prove, for example with an exact computation
in the Brownian setting , that
the continuous critical curve is sensitive to a depinning term, i.e.,
$K_c^{\beta}<K_c^{0}$ for certain $\beta<0$, then Theorem \ref{theo3} will entail that the same
$\beta<0$ satisfies
$h_c^{\lambda \beta}(\lambda)=K_{c}^{\beta} \lambda (1+o(1))$. This would prove that
$\mathcal{L}$ shrinks under the influence of a depinning term, at least for $\lambda$ small.

\subsubsection{The homopolymer with adsorption}
By fixing $\lambda=1$ and $w_{i}\equiv 0$ for $i\geq 1$, we can model a homopolymer
instead of a copolymer. Effectively, in this case the polymer only consists of hydrophobic monomers, and
its Hamiltonian is given by
\begin{equation}\label{Hamil}
h \sum_{i=1}^{N}
\Lambda_{i}+ \beta \sum_{j=-K}^{K} \sum_{i=1}^{N} \gamma_{i}^{j} \ind_{\{S_{i}=j\}}.
\end{equation}
%In what follows, we will refer to this model as the $h$-model. Apart from the fact
%that this $h$-model corresponds to a different physical situation (homopolymer instead of copolymer), there are
%two principal reasons to study it. First,
This type of model, which we call $h$-model, with a pinning term at the interface
in competition with a repulsion effect (given here by $h \sum_{i=1}^{N}\Lambda_{i}$), has already been
investigated
in the literature (see \cite{JSW}, or \cite{Der}). It has been proven, for instance,
that some properties of the $h$-model can be extended to the wetting model by letting the parameter $h$
tend to $\infty$ (see \cite{Petr2}).

The free energy of the $h$-model is
denoted by $\Phi(\beta,h)$ and the localization condition is as before: $(\beta,h)\in \mathcal{L}$ when
$\Phi(\beta,h)>h$ and $(\beta,h)\in \mathcal{D}$ when
$\Phi(\beta,h)=h$. The critical curve of the $h$-model, which separates the $(h,\beta)$-plane
into a localized and a delocalized phase, is denoted by $\kappa_{c}(\beta)$.
%covers the   In the wetting model for instance, the repulsion is given
%by the fact that the involved random walk is conditioned to stay positive. We will prove in chapter $4$
%that this model can be seen as the limit of the $h$-model when $h$ tends to infinity. In that way, we
%will translate
%some results concerning the $h$-model to the wetting one. We show also in Chapter $4$ that
%this $h$-model
This curve is increasing, convex and satisfies $\kappa_{c}(0)=0$.

At this stage, we must recall that $w$ is assumed to satisfy $\mathbb{E}(w_1^2)=1$ in Theorem \ref{theo3}.
Therefore, this theorem can not be applied directly to the $h$-model. However, the proof of
Theorem \ref{theo3}, that we give in Section \ref{sec3}, can easily be extended to the $h$-model, so
that \eqref{theo31} can be restated in this case as
\begin{equation}
\lim_{a\to \infty} \frac{1}{a^2}\Phi(a\beta,a^2h)=\tilde{\Phi}(\beta\Sigma,h),
\end{equation}
where $\tilde{\Phi}(\beta\Sigma,h)$ denotes the free energy of the continuous limit of the $h$-model.
The hamiltonian
of this continuous limit is given by
\begin{equation}\label{Ham2}
h \int_{0}^{t} \Lambda_{s} ds +\beta \Sigma L_{t},
\end{equation}
which is remarkable because the disorder disappears. Thus,
we can compute explicitly some quantities related to $\tilde{\Phi}$. For instance, we state the following
proposition for the case $\Sigma=1$.
\begin{proposition}\label{proph}
Let $\beta\in \mathbb{R}$ and $h\geq 0$. Then,
$$\tilde{\Phi}(\beta,h)=h\ \ \ \text{if}\ \ \ h\geq \beta^{2}\ \ \ \ \text{and}\ \ \ \ \tilde{\Phi}
(\beta,h)
=\frac{h^{2}}{2 \beta^{2}}+ \frac{\beta^{2}}{2}\ \ \text{if}\ \ \ h<\beta^{2}.$$
\end{proposition}
%The localization condition remains the same, i.e.
%$(\beta,h)\in \mathcal{L}$ when $\Phi(\beta,h)>h$ or $\tilde{\Phi}(\beta,h)>h$ in the continuous case.
\noindent Since
$h^{2}/(2 \beta^{2})+ \beta^{2}/2>h$ when $h<\beta^{2}$, we obtain the continuous critical curve, i.e.,
$\tilde{\kappa}_{c}(\beta)=\beta^{2}$ for $\Sigma=1$ (see Fig. $2$).

\vspace{1.5cm}

\setlength{\unitlength}{0.23cm}
\begin{picture}(20,15)(-22,0)

\put(-13,0){\vector(1,0){40}}
\put(0,0){\vector(0,1){15}}

\put(-.8,-1.3){$0$}

\put(32,-.3){$\beta$}

\put(-1,16){$h$}
{\thicklines
   \qbezier(0,0)(11,0)(15.6,12.5)}

% la legende
\put(-18,17){\small
             Fig.\ 2:}%{\thicklines
\put(-5.3,9){$\mathcal{D}:$}
\put(-5.3,6.5){$\tilde{\Phi}(\beta,h)=h$}
\put(18,9){$\mathcal{L}:$}
\put(18,6.5){$\tilde{\Phi}(\beta,h)=\frac{h^{2}}{2\beta^{2}}+\frac{\beta^{2}}{2}$}
\put(11.5,14){$\tilde{\kappa}_{c}(\beta)=\beta^{2}$}
\end{picture}

$\phantom{aaa}$\\

Thanks to Proposition \ref{proph} we can give the asymptotic behavior, as $\beta$ tends to $0$,
of some quantities linked to the discrete model.
For instance, for the general $h$-model, i.e. with $\Sigma$ not necessarily equal to $1$,
we can state the equivalent of \eqref{corol}, that is
\begin{equation}\label{rebab}
\lim_{\beta \to 0}\frac{\kappa_{c}(\beta)}{\beta^{2}}=\Sigma^{2}
.\end{equation}
Proving \eqref{rebab} requires to restate Theorem \ref{theo} (introduced below) for the $h$-model. This
does not present any further difficulty, that is why we will not give the details here.
Notice that the limit \ref{rebab} is conform to our intuition that a stronger pinning along the interface
enlarges the localized area and, consequently,
increases the curvature of the critical curve at the origin.
It is also confirmed by the bounds on the critical curve found in \cite{Petr2}.% the Chapter $4$ (in
%that particular case $K$ is equal to $1$, $\gamma$ satisfies $E(\gamma_{1})=1$ and
%$h_{c}(\beta)=(1+o(\beta)) \beta^{2}$
%as $\beta$ tends to $0$).

Still with Proposition \ref{proph}, we can differentiate $\tilde{\Phi}(h,\beta)$ with respect to $\beta$ and we find
the asymptotic behavior of the reward average in the weak coupling limit. Indeed, if $h<\beta^{2}$, then
by convexity of $\Phi_{N}$ in $\beta$ we can state that, a.s. in $\gamma$,
\begin{equation*}
\lim_{a\to 0}\,\lim_{N \to \infty}\,\textstyle{\frac{1}{a N}E_{N,a\beta}^{a^{2}h,w}
\Big[\sum_{j=-K}^{K}\sum_{i=1}^{N}\gamma_{i}^{j}
\ind_{\{S_{i}=j\}}\Big]=\beta-\frac{h^{2}}{\beta^{3}}}.
\end{equation*}
The same derivative with respect to $h$ gives an approximation, for $a$ small, of the time proportion
spent by the polymer under the interface, i.e.,
\begin{equation*}
\lim_{a\to 0}\,\lim_{N \to \infty}\,\textstyle{E_{N,a\beta}^{a^{2}h,w}\Big[\frac{\sum_{i=1}^{N} \Delta_i}{N}\Big]
=\frac{\beta^{2}-h}{2\beta^{2}}}.
\end{equation*}

\subsection{Organization of the paper}

In section \ref{sec2}, we will state and prove some technical results that turn out to be useful in the
proof of Theorem \ref{theo3}. More precisely, in section \ref{sec21}, we consider the local time spent by the
random walk in a finite layer around the interface after $N$ steps. We rescale the later by $\sqrt{N}$ and
we prove its convergence, in terms of exponential moments, towards the local time spent at the origin
by the Brownian motion between times $0$ and $1$.
%In section \ref{sec22} we introduce the so called excess free energy, which allows to
In section \ref{sec23} we introduce
the Theorem \ref{theo}, from which Theorem \ref{theo3} will be deduced. Theorem \ref{theo} is essentially
technical and
consists in comparing the continuous
free energy and the discrete free energy when the coupling is weak. Finally, in Section \ref{finir}, we provide
some conditions of finiteness for $h_c$ and $\tilde{h_c}$.

Section \ref{sec3} is essentially dedicated to the proof of Theorems \ref{theo3} and \ref{theo}. Thus,
in section \ref{sec32} we explain how Theorem \ref{theo3} is deduced from Theorem \ref{theo}, whereas
the rest of Section $3$ is dedicated to the proof of Theorem \ref{theo}.

Section \ref{sec4} is an appendix dedicated to the exact computation of $\tilde{\Phi}$ asserted in Proposition
\ref{proph}.

\section{Preparation}\label{sec2}

\subsection{Technical Lemma}\label{sec21}

\begin{lemma}\label{lemsi}
For every $K \in \mathbb{N}$ and every $(f_{-K},f_{-K+1},\dots,f_{K})$ in $\mathbb{R}^{2K+1}$
the following convergence occurs:
\begin{equation}
\lim_{N\to \infty}E\bigg[\exp\bigg(\frac{1}{\sqrt{N}}\sum_{j=-K}^{K}f_{j}
\sum_{i=1}^{N}\ind_{\{S_{i}=j\}}
\bigg)\bigg]=E\bigg[\exp\bigg(\bigg(\sum_{j=-K}^{K}f_{j}\bigg)\ \  L_{1}^{0}\bigg)\bigg]
,\end{equation}
where $L_{1}^{0}$ is the local time in $0$ of a Brownian motion $\left(B_{s}\right)_{s\geq 0}$ between $0$ and $1$.
\end{lemma}
\begin{proof}
First, we prove the following intermediate result. For every $K \in \mathbb{N}$
\begin{equation}\label{conv}
\lim_{N\to \infty}\textstyle{\frac{1}{\sqrt{N}}}\sum_{j=-K}^{K} f(j) \sum_{i=1}^{N}\ind_{\{S_{i}=j\}}
\stackrel{\text{Law}}{=}
\big(\sum_{j=-K}^{K} f_{j}\big)\ \ L_{1}^{0}
.\end{equation}
For simplicity, we only prove that
$\frac{1}{\sqrt{N}}\big(\sum_{i=1}^{N}\ind_{\{S_{i}=0\}},\sum_{i=1}^{N}\ind_{\{S_{i}=1\}}\big)$
converges in law to $\left(L_{1}^{0},L_{1}^{0}\right)$ as $N \uparrow\infty$.
The proof for $2K+1$ levels is exactly the same.
For this convergence in law, we use a result of \cite{Rev}, saying that
we can build, on the same probability space $(W,\mathcal{A},P)$, a simple random walk
$\left(S_{i}\right)_{i\geq 0}$ and
a Brownian motion $(B_{s})_{s\geq 0}$ such that $P$ almost surely
\begin{equation}\label{csa}
\lim_{n\to \infty}  \sup_{j\in\{0,1\}}\textstyle{\frac{1}{\sqrt{n}}}\big{|}U_{n}^{j}-L_{n}^{j}\big{|}=0
\end{equation}
with $U_{n}^{j}=\sum_{i=1}^{n}\ind_{\{S_{i}=j\}}$ and $L_{n}^{x}$ the local time in $x$ of $B$
between $0$ and $n$.
The equation \eqref{csa} implies that $\frac{1}{\sqrt{n}}\ (U_{n}^{0}-L_{n}^{0})$ and
$\frac{1}{\sqrt{n}}\ (U_{n}^{1}-L_{n}^{1})$
tend a.s. to $0$ as $n \uparrow\infty$. Therefore, the proof of \eqref{conv} will be completed
if we show that $\frac{1}{\sqrt{n}}\ (L_{n}^{0},L_{n}^{1})$
converges in law to $(L_{1}^{0},L_{1}^{0})$. By the scaling property of Brownian motion, we obtain that,
for every $n\geq 1$, $\frac{1}{\sqrt{n}}\ (L_{n}^{0},L_{n}^{1})$
has the same law as $(L_{1}^{0},L_{1}^{1/\sqrt{n}})$. Thus, since $L_{1}^{x}$ is a.s. continuous in $x=0$,
we obtain
immediately the a.s. convergence of $(L_{1}^{0},L_{1}^{1/\sqrt{n}})$ towards $(L_{1}^{0},L_{1}^{0})$.
This a.s. convergence implies the convergence in law and \eqref{conv} is proven.

Since the function $\exp(x)$ is continuous, \eqref{conv} gives us the convergence
in law of $W_{N}=\exp\big(\frac{1}{\sqrt{N}}\sum_{j=-K}^{K}f_{j}\sum_{i=1}^{N}\ind_{\{S_{i}=j\}}\big)$
to $\exp\big(\big(\sum_{j=-K}^{K} f_{j}\big)L_{1}^{0}\big)$ as $N \uparrow\infty$.
The uniform integrability of the sequence $\big(W_{N}\big)_{N\geq 1}$ will therefore be sufficient to complete
the proof of Lemma \ref{lemsi}.

We will obtain this uniform integrability if we can prove that $\sup_{N\geq 1}E(W_{N}^2)<\infty$.
By the H\"older inequality, it is sufficient to prove that for every $b>0$ and every $j\in\mathbb{Z}$
we have the inequality
\begin{equation}\label{kn}
\textstyle{\sup_{N\geq 1}E\Big(\exp\big(\frac{b}{\sqrt{N}}\sum_{i=1}^{N}\ind_{\{S_{i}=j\}}\big)\Big)<\infty}.
\end{equation}
We let $k_N=\sum_{i=1}^{N}\ind_{\{S_{i}=0\}}$ and $\tau_j=\inf\{n\geq 1 : S_n=j\}$. Thus by the
Markov property we can write
\begin{equation}\label{kn2}
\textstyle{E\Big(\exp\big(\frac{b}{\sqrt{N}}\sum_{i=1}^{N}\ind_{\{S_{i}=j\}}\big)\Big)
\leq E\Big(\exp\big(\frac{b }{\sqrt{N}}\,\ind_{\{\tau_j \leq N\}}\sum_{i=\tau_j}^{N}\ind_{\{S_{i}=j\}}\big)\Big)}
\leq E\big(e^{\frac{b}{\sqrt{N}}(1+k_N)}\big),
\end{equation}
and it just remains to prove that for every $b>0$
the sequence $\big(E\big[\exp(b k_{N}/\sqrt{N})\big]\big)
_{N\geq 0}$ is bounded from above independently of $N$. To that aim, we notice that $k_{N}\leq k_{2N}\leq N$ and
 write the obvious inequality
\begin{equation}\label{eq1}
\textstyle{E\big[\exp(b k_{2N}/\sqrt{N})\big]\leq \sum_{k=0}^{\lfloor\sqrt{N/2}\rfloor}\, e^{\sqrt{2}\, b (k+1)}\,
P\big(k_{2N} \in \big[k \sqrt{2N},(k+1)
\sqrt{2N}\big[\big)}.
\end{equation}
With the help of \cite{Fel} we can compute an upper bound of $P\big(k_{2N} \in \big[k \sqrt{2N},(k+1)
\sqrt{2N}\big[\big)$. Indeed, for every $k \leq \big\lfloor\sqrt{N/2}\big\rfloor$ we obtain
\begin{equation}\label{eq2}
P\big(k_{2N} \in \big[k \sqrt{2N},(k+1)
\sqrt{2N}\big[\big)
\leq \sum_{j=\lfloor k \sqrt{2N}\rfloor}^{\max(\,\lfloor(k+1)\sqrt{2N}\rfloor,N)}
P(S_{2N}=0)
\textstyle{\frac{(1-\frac{1}{N})\dots(1-\frac{j-1}{N})}
{(1-\frac{1}{2N})\dots(1-\frac{j-1}{2N})}}
.\end{equation}
The function $x\to \log(1-x)+x$  is decreasing on $[0,1)$ and consequently, for every\  $j\in
\{\lfloor k \sqrt{2N}\rfloor,\dots,\max(\,\lfloor(k+1)\sqrt{2N}\rfloor,N)\}$, we have
$\log(1-j/N)-\log(1-j/2N)\leq -j/2N$. Therefore,
\begin{equation*}
\textstyle{\frac{\left(1-\frac{1}{N}\right)\dots\left(1-\frac{j-1}{N}\right)}
{\left(1-\frac{1}{2N}\right)\dots\left(1-\frac{j-1}{2N}\right)}\leq \exp
\big(\sum_{i=1}^{j-1}-\frac{i}{2N}\big)=\exp\big(-\frac{j(j-1)}
{4N}\big)\leq \exp\big(-\frac{(k-1)^{2}}{2}
\big)}.
\end{equation*}\\
Moreover $\big\lfloor(k+1)\sqrt{2N}\big\rfloor-\big\lfloor k \sqrt{2N}\big\rfloor\leq \sqrt{2N}+1$ and
there exists a constant $c>0$ such that $P\left(S_{2N}=0\right)\leq c/\sqrt{2N}$
for every $N\geq 1$ . That is why, the equation \eqref{eq2} becomes
$$P\big(k_{2N} \in \big[k \sqrt{2N},(k+1)
\sqrt{2N}\big]\big)\leq 2c \exp\left(-(k-1)^{2}/2\right).$$
This results allows us to rewrite \eqref{eq1} as
\begin{equation*}
\textstyle{E}\big[\exp(b k_{2N}/\sqrt{N})\big]\leq \sum_{k=0}^{\infty} 2\,c\, e^{b (k+1)}\, e^{-\frac{(k-1)^{2}}{2}}
,\end{equation*}
and the r.h.s. of this inequality is the sum of a convergent series. Therefore, the sequence
$(W_N)_{N\geq 0}$ is uniformly integrable and the proof of Lemma \ref{lemsi} is completed.%\qed
\end{proof}

\subsection{Excess free energies} \label{sec23}
%\begin{rem}
We define the quantities
$\Psi_{N}(\beta, \lambda, h)=\Phi_{N}(\beta, \lambda, h)-\lambda h$ and
$\tilde{\Psi}_{t}(\beta, \lambda, h)=\tilde{\Phi}_{t}(\beta, \lambda, h)-\lambda h$. They
converge respectively to $\Psi(\beta, \lambda, h)
=\Phi(\beta, \lambda, h)-\lambda h$ and $\tilde{\Psi}(\beta, \lambda, h)=
\tilde{\Phi}(\beta, \lambda, h)-\lambda h$, which are called {\it excess free energies} of the polymer.
Therefore, to decide whether the polymer is localized or not,
it suffices to compare $\Psi$ or $\tilde{\Psi}$ with $0$.
Moreover, since  $\sum_{i=1}^{N} (w_i+h)=hN+o(N)$ when $N\uparrow \infty$, $\mathbb{P}$-a.s.,
we can subtract this quantity
from the Hamiltonian \eqref{Ham1} and associate $\Psi_N$ with
\begin{equation*}
H_{N,\beta,\lambda,h}^{w,\gamma}=-2 \lambda \sum_{i=1}^{N} (w_{i}+h) \Delta_{i}+
\beta \sum_{j=-K}^{K}\sum_{i=1}^{N} \gamma_{i}^{j}\ \ind_{\{S_{i}=j\}},
\end{equation*}
with $\Delta_{i}=1$ if $\Lambda_{i}=-1$ and  $\Delta_{i}=0$ otherwise.
Similarly, $\tilde{\Psi}_{t}(\beta, \lambda, h)$
is associated with
\begin{equation*}
\tilde{H}_{t,\beta}^{R}=-2 \lambda \int_{0}^{t} \ind_{\{B_{s}<0\}} (dR_{s}+h ds)
+\beta L_{t}^{0},
\end{equation*}
and $\Psi$ and
$\tilde{\Psi}$ are continuous, separately convex and
non-increasing
in $h$. Moreover $\tilde{\Psi}$ is non-decreasing in $\beta$.

\subsection{Technical Theorem}\label{sec23}
\begin{rem}\rm{Stating Theorem \ref{theo} requires a slight modification of the Hamiltonian. In fact,
let $(\beta_{1},\beta_{2})\in\mathbb{R}^2$ and define
\begin{align*}
I_{1}=\{j\in\{-K,\dots,K\}\colon\, \mathbb{E}(\gamma_{1}^{j})>0\}\ \ \text{and}\ \
I_{2}=\{j\in\{-K,\dots,K\}\colon\, \mathbb{E}(\gamma_{1}^{j})<0\}.
\end{align*}
Then, if $\mathbb{E}(\gamma_{1}^{j})\neq 0$ for every $j\in\{-K,\dots,K\}$, we define
\begin{equation}
H_{N,\beta_{1},\beta_{2},\lambda,h}^{w,\gamma}=\beta_{1}\sum_{j\in I_{1}}\sum_{i=1}^{N}
\gamma_{i}^{j} \ind_{\{S_{i}=j\}}+\beta_{2}\sum_{j\in I_{2}}\sum_{i=1}^{N}
\gamma_{i}^{j} \ind_{\{S_{i}=j\}}+\lambda \sum_{i=1}^{N}(w_{i}+h) \Lambda_{i}
.\end{equation}
The associated free energy $\Psi(\beta_{1},\beta_{2},\lambda,h)$ is defined as in Proposition \ref{prop1}, and satisfies
$\Psi(\beta,\lambda,h)=\Psi(\beta,\beta,\lambda,h)$. Thus, in what follows, we will
use the notation $\Psi(\beta_{1},\beta_{2},\lambda,h)$ if $\beta_{1}\neq \beta_{2}$, otherwise
we will use $\Psi(\beta,\lambda,h)$. We let $\Sigma=\Sigma_{1}+\Sigma_{2}$, with
$\Sigma_{1}= \sum_{j\in I_{1}}\mathbb{E}
(\gamma_{1}^{j})$ and $\Sigma_{2}=\sum_{j\in I_{2}}\mathbb{E}
(\gamma_{1}^{j})$}.
\end{rem}
\begin{theorem}\label{theo}
%If $\mathbb{E}\big(\gamma_{1}^{j}\big)>0$ for every $j\in \{-K,\dots,K\}$,
Suppose $\mathbb{E}(\gamma_{1}^{j})\neq 0$ for every $j\in\{-K,\dots,K\}$.
If $\beta_{1}\neq 0$, $\beta_{2}\neq 0$, and $(\mu_{1},\mu_{2}) \in \mathbb{R}^{2}$ satisfy
$$\mu_{1}>\beta_{1}\Sigma_{1}+\beta_{2}\Sigma_{2}>\mu_{2},$$
and $\rho>0$, $h>0$, $h^{'}\geq 0$,
$\lambda>0$
satisfy $(1+\rho)h^{'}<h$, then there exists $a_{0}>0$ such that for every $a<a_{0}$
\begin{align}\label{eqgen}
\frac{1}{a^{2}}\Psi\left(a\beta_{1}, a\beta_{2}, a\lambda , a h\right)&\leq (1+\rho)\  \tilde{\Psi}
\big(\mu_{1}, \lambda, h^{'}\big)\\
\nonumber
\tilde{\Psi}(\mu_{2}, \lambda, h)&\leq \frac{1+\rho}{a^{2}}\  \Psi\big(a
\beta_{1}, a\beta_{2}, a \lambda, a h^{'}\big)
.\end{align}
\end{theorem}

\subsection{Conditions of finiteness for $h_c$ and $\tilde{h}_c$}\label{finir}

\subsubsection{One particular discrete model}

We give here some
more details about one particular case, namely $K=0$ and $\gamma_1^0=1$.
We let $\zeta(\lambda)=\log(\mathbb{E}(\exp(\lambda w_1))$. The Jensen's inequality allows us to write
\begin{equation}\label{aaa}
\textstyle \Psi(\beta,\lambda,h)\leq
\lim_{N\to\infty} \frac{1}{N} \log E\Big[\exp\big(\sum_{i=1}^{N} (\zeta(-2\lambda)-2\lambda h) \Delta_{i}+
\beta \sum_{i=1}^{N} \ind_{\{S_{i}=0\}}\big)\Big].
\end{equation}
The limit in the right-hand side of \eqref{aaa} is computed in \cite{Petr} and is equal to $0$ for $h$
large enough, as long as $\beta<\log 2$. This means that $h_c^{\beta}(\lambda)<\infty$ for $\beta<\log 2$.
In the same spirit we can let $h$ tend to $\infty$ and write the lower bound
\begin{equation}\label{bbb}
\textstyle \Psi(\beta,\lambda,h)\geq
\lim_{N\to\infty} \frac{1}{N} \log E\Big[\exp\big(
\beta \sum_{i=1}^{N} \ind_{\{S_{i}=0\}}\big)\ \ind_{\{S_i\geq 0\  \forall i\in\{1,\dots,N\}\}}\Big].
\end{equation}
The r.h.s. of \eqref{bbb} is strictly positive for $\beta>\log 2$, and therefore $h_c^{\beta}(\lambda)=\infty$
when $\beta>\log 2$.

\subsubsection{The continuous case}

In the continuous case we can assert the following general result.
\begin{proposition}\label{propfin}
For every $\beta \in \mathbb{R}$ we have $K_c^{\beta}<\infty$. As a consequence,
for every $\beta\in\mathbb{R}$ and $\lambda>0$ we have $\tilde{h}_{c}^{\beta}(\lambda)<\infty$.
\end{proposition}
The proof of this proposition involves the discrete case mentioned above and
the Theorem \ref{theo}.

\vspace {0.5 cm}

\section{Proof of theorems and propositions}\label{sec3}

\subsection{Proof of Proposition \ref{propfin}}
In this section we assume that Theorem \ref{theo} is satisfied. Since $\beta \to K_c^\beta$ is non decreasing in
$\beta$, the proof of Proposition \ref{propfin} will be completed if we can show
that $K_c^\beta<\infty$ for all $\beta>0$. Therefore, we let $\beta>0$ and for any $h\geq 0$, we
let
\begin{equation}
\textstyle{\Upsilon(\beta,h)=\lim_{N\to\infty} \frac{1}{N} \log E\Big[\exp\big(-2h\sum_{i=1}^{N} \Delta_{i}+
\beta \sum_{i=1}^{N} \ind_{\{S_{i}=0\}}\big)\Big]}.
\end{equation}
It is proven in \cite{Petr} that $\Upsilon(\beta,h)>0$ when $h<\kappa_c(\beta)$ and $\Upsilon(\beta,h)=0$
when $h\geq \kappa_c(\beta)$. The critical value is also computed in \cite{Petr}, i.e.,
\begin{equation}\label{cc1}
\kappa_c(\beta)=-(1/4)\log(1-4(1-\exp(-\beta))^2).
\end{equation}
We recall the particular discrete case introduced in section \ref{dud}, namely $K=0$ and $\gamma_1^0=1$.
We assume also that $w_1$ is a Bernouilli trial taking the values $1$ and $-1$ with probability $1/2$.
We let $\beta>0$, $h>0$ and we can apply the second inequality
of Theorem \ref{theo} to this particular discrete model with the parameters $\rho=1/2$, $\mu_2=\beta$,
$\beta_1=\beta_2=2\beta$ and $h'=h/2$. Since $\Sigma_1=1$ and $\Sigma_2=0$ in this case
we obtain for $a$ small enough
\begin{equation}\label{duda}
\tilde{\Psi}(\beta,1, h)\leq \frac{1+\frac{1}{2}}{a^{2}}\  \Psi\big(2a
\beta,a, a h/2\big).
\end{equation}
Moreover, the equation \eqref{aaa} gives $\Psi(2a\beta,a, ah/2)\leq
\Upsilon(2a\beta,-\zeta(-2a)/2+a^2 h/2)$. The equation \eqref{cc1}
gives $\kappa_c(2a\beta)=4a^2\beta^2+o(a^2)$ whereas $-\zeta(-2a)/2+a^2 h/2=a^2(h-2)/2+o(a^2)$. Therefore,
by choosing $h$ large enough and $a$ small enough we have that the r.h.s. of \eqref{duda} is equal to $0$.
This shows that $K_c^\beta<\infty$ for all $\beta>0$
and
the proof of Proposition \ref{propfin} is
completed.

\subsection{Proof of Theorem \ref{theo3}} \label{sec32}
\subsubsection{Part 1: proof of the convergence \eqref{theo31}}
In this section, we prove that the convergence \eqref{theo31} is a consequence of Theorem \ref{theo}.
%Recall that in the case $\beta=0$ the Theorem is already proven in \cite{BDH}.
%Therefore we consider $\beta\neq 0$.
This proof is divided into $3$ steps.
In the first step, we show that \eqref{theo31} is satisfied when $\lambda>0$, $h>0$, $\beta\neq 0$
and every pinning reward
$\gamma_{1}^{j}$ has a non zero average. In the second step, we prove that the result can be extended to the
case in which some $\gamma_{1}^{j}$ have a zero average and consequently to the case $\beta=0$. Finally,
in the last step, we will consider the case $h=0$. We recall that proving \eqref{theo31} with $\Phi$
and $\tilde{\Phi}$ or $\Psi$ and $\tilde{\Psi}$ is completely equivalent.

\subsubsection{\bf{Step} $1$}

First, we consider the case $\lambda>0$, $h>0$, $\beta\neq 0$ and
$\mathbb{E}(\gamma_{1}^{j})\neq 0$ for every $j\in \{-K,\dots,K\}$. We can apply
the first inequality of Theorem \ref{theo} with the parameters
$\rho=1/n$, $h'=h/(1+1/n)^{2}$, $\beta_{1}=\beta_{2}=\beta$ and  $\mu_{1}(v)=\beta \Sigma_{1}+\beta
\Sigma_{2}+1/v$ ($n$ and $v\in \mathbb{N}-\{0\}$).
It gives, for every integers $n$ and $v$ strictly positive, that
\begin{equation}\label{V_{1}}
\textstyle{\limsup_{a \to 0}}\, \frac{1}{a^{2}}\, \Psi\big(a \beta, a \lambda, a h\big)\leq
(1+\frac{1}{n}) \tilde{\Psi}\big(\mu_{1}(v), \lambda,
\frac{h}{(1+1/n)^{2}}\big)
.\end{equation}
At this stage, we let successively $n$ and $v$ tend to $\infty$, and, by continuity of $\tilde{\Psi}$ in $h$
and $\beta$ we obtain $\limsup_{a \to 0} 1/a^{2} \Psi\left(a \beta, a \lambda, a h\right)\leq
\tilde{\Psi}(\beta \Sigma, \lambda, h)$. The lower bound is proven with the second inequality of
Theorem \ref{theo}.
Indeed, if we choose $\mu_{2}(v)=\beta \Sigma_{1}+\beta \Sigma_{2}-1/v$ and keep the other notations,
we obtain
\begin{equation}\label{V_{2}}
\textstyle{\tilde{\Psi}}\big(\mu_{2}(v), \lambda, h(1+\frac{1}{n})^{2}\big)\leq
(1+\frac{1}{n}) \liminf_{a\to 0}\,\frac{1}{a^{2}} \Psi(a \beta, a \lambda, a h).
\end{equation}
We let $n \uparrow\infty$, and after, we let $v\uparrow\infty$. In that way, we can conclude
that $\lim_{a \to 0} 1/a^{2}\Psi\left(a \beta, a \lambda, a h\right)= \tilde{\Psi}(\beta \Sigma, \lambda, h)$
which implies \eqref{theo31}.

\subsubsection{\bf{Step} $2$}

We prove the convergence \eqref{theo31} when $\lambda>0$, $h>0$, $\beta\neq 0$ and there exists $j\in \{-K,\dots,K\}$
such that $\mathbb{E}(\gamma_{1}^{j})=0$. For that, we choose $\mu>0$ and small enough, such that,
$\mathbb{E}(\gamma_{i}^{j}+\mu)\neq 0$ for every $j\in \{-K,\dots,K\}$. With these new variables we
can use the result of Step $1$ with $\Sigma_{\mu}=\Sigma+(2K+1)\mu$.
Since
the free energy $\Psi_{\mu}$ associated with the variables $\gamma_{i}^{j}+\mu$ is larger than $\Psi$, we obtain
\begin{equation*}
\limsup_{a\to 0} \frac{1}{a^{2}}\Psi(a\beta,a \lambda,a h)\leq
\lim_{a\to 0} \frac{1}{a^{2}}\Psi_{\mu}(a\beta,a \lambda,a h)=\tilde{\Psi}(\beta (\Sigma+(2K+1)\mu), \lambda,h).
\end{equation*}
As $\tilde{\Psi}$ is continuous in $\beta$, we let $\mu \downarrow 0$ and write
$\limsup_{a\to 0} 1/a^{2}\Psi(a\beta, a \lambda, a h)\leq \tilde{\Psi}(\beta \Sigma, \lambda, h)$.
Thus, it suffices to do the same computation with $-\mu<0$, and we obtain the other inequality, i.e.,
\begin{equation*}
\liminf_{a\to 0} \frac{1}{a^{2}}\Psi(a\beta,a \lambda,a h)\geq
\lim_{\mu \to 0}\tilde{\Psi}(\beta (\Sigma-(2K+1)\mu), \lambda,h)=\tilde{\Psi}(\beta \Sigma, \lambda, h).
\end{equation*}
Therefore, we can say that
$\lim_{a\to 0} \frac{1}{a^{2}}\Psi(a\beta,a \lambda,a h)=
\tilde{\Psi}(\beta \Sigma, \lambda, h)$.

As a consequence, \eqref{theo31} is satisfied when the variables $\gamma_i^j$
are all equal to $0$. Therefore, it is also satisfied when $\beta=0$.

%At this stage, we denote $S=\{j\in \{-K,\dots,K\}\,:\, \mathbb{E}(\gamma_{1}^{j})=0\}$ and we denote by
%$H'(S)$ the Hamiltonian
%\begin{equation*}
%\lambda \sum_{i=1}^{t/a^{2}} (w_{i}+h)
%\Delta_{i}+ \beta \sum_{j\in\{-K,\dots,K\}-S}\ \sum_{i=1}^{t/a^{2}}\gamma_{i}^{j}\ \ind{\{S_{i}=j\}}+
%\beta_{1} \sum_{j\in S}\sum_{i=1}^{t/a^{2}}\gamma_{i}^{j}\ \ind{\{S_{i}=j\}}
%,\end{equation*}
%and by $\Phi^{'}(\beta,\beta_{1},\lambda,h)$ the associated free energy. By independence of the $\gamma_{i}^{j}$
%we obtain easily $\partial_{\beta_{1}}\Phi'(\beta,0,\lambda,h)=0$ and $\Phi'$ is convex in $\beta_{1}$,
%therefore $\Phi'$ is not
%decreasing in $\beta_{1}$. That is why we can write
%\begin{equation}\label{tera}
%\liminf_{a\to 0} \frac{1}{a^{2}}\Phi(a\beta, a \lambda, a h)=\liminf_{a\to 0}
%\frac{1}{a^{2}}\Phi'(a\beta,a\beta,a \lambda,a h)
%\geq \liminf_{a\to 0} \frac{1}{a^{2}}\Phi'(a\beta,0, a \lambda,a h)
%,\end{equation}
%and we notice that theorem $3$ with strictly positive parameters, can be applied to the rhs of \eqref{tera}
%because the $\gamma_{i}^{j}$ of $0$
%average (with $j\in S$) do not play any role in the computation of $\Phi'(a\beta,0,a \lambda,a h)$. So, we obtain
%$\liminf_{a \to 0} 1/a^{2}\Phi(a\beta,a \lambda,a h)\geq \tilde{\Phi}(\beta \Sigma, \lambda,h)$
%and it completes the proof of theorem $3$ in this case.

\subsubsection{\bf{Step} $3$}

It remains to prove \eqref{theo31} when $h=0$. Since $\Psi$ and $\tilde{\Psi}$
are non increasing in $h$, \eqref{theo31} with $\lambda>0$, $h>0$ and $\beta\in\mathbb{R}$
(proven in Step $2$) implies
\begin{equation*}
\liminf_{a \to 0}\frac{1}{a^{2}}\Psi(a\beta,a \lambda,0)
\geq \liminf_{a \to 0}\frac{1}{a^{2}}\Psi(a\beta,a \lambda,ah)= \tilde{\Psi}(\beta \Sigma, \lambda,h).
\end{equation*}
We let $h\downarrow 0$ and by continuity of $\tilde{\Psi}$ in $h$ we obtain
\begin{equation*}
\liminf_{a \to 0}\frac{1}{a^{2}}\Phi(a\beta,a \lambda,0)
=\liminf_{a \to 0}\frac{1}{a^{2}}\Psi(a\beta,a \lambda,0)\geq \tilde{\Psi}(\beta \Sigma, \lambda,0)
=\tilde{\Phi}(\beta \Sigma, \lambda,0).
\end{equation*}
To prove the opposite inequality, we just notice that $\Phi$ is non decreasing in $h$. Effectively
%\begin{equation}\label{deriv}
%\textstyle{\frac{\partial \Phi_N}{\partial h}}\big|_{(\beta,\lambda,0)}=
%\mathbb{E}\Big[\sum_{(\Lambda_{1},\dots,\Lambda_{N})\in \{-1,1\}^{N}}
%P(\Lambda_{1},\dots,\Lambda_{N})\ \frac{\lambda \sum_{i=1}^{N}\Lambda_{i}}{N}\  \exp\big
%(\lambda \sum_{i=1}^{N}w_{i}
%\Lambda_{i}+\beta \dots\big)\Big],
%\end{equation}
\begin{equation}\label{deriv}
\textstyle{\frac{\partial \Phi_N}{\partial h}}\big|_{(\beta,\lambda,0)}=\mathbb{E}
\Big[E_{N,\beta,\lambda,h}^{w,\gamma}\big(\frac{\lambda}{N}\sum_{i=1}^N \Lambda_i\big)\Big],
\end{equation}
and by symmetry of the laws of the random walk and of the variables $\{w_{i}\}_{i=1,2,\dots}$, we can transform
$w_{i}$ in $-w_{i}$,
and $(\Lambda_{1},\dots,\Lambda_{N})$ in $(-\Lambda_{1},\dots,-\Lambda_{N})$, without changing \eqref{deriv}.
It gives
\begin{align*}
\textstyle{\frac{\partial \Phi_N}{\partial h}}\big|_{(\beta,\lambda,0)}=
-\frac{\partial \Phi_N}{\partial h}\Big|_{(\beta,\lambda,0)}.%\mathbb{E}
%\Big[\sum_{\Lambda_{1},\dots,\Lambda_{N}}
%P(-\Lambda_{1},\dots,-\Lambda_{N})\ \frac{-\lambda \sum_{i=1}^{N}\Lambda_{i}}{N}\  \exp\big
%(-\lambda \sum_{i=1}^{N}-w_{i}
%\Lambda_{i}+\beta \dots\big)\Big]
\end{align*}
%\begin{align*}
%\hspace{0.65cm}\textstyle{=}\,\mathbb{E}\Big[\sum_{\Lambda_{1},\dots,\Lambda_{N}}
%P(\Lambda_{1},\dots,\Lambda_{N})\ \frac{-\lambda \sum_{i=1}^{N}\Lambda_{i}}{N}\  \exp\big
%(\lambda \sum_{i=1}^{N}w_{i}
%\Lambda_{i}+\beta \dots\big)\Big]=-\frac{\partial \Phi_N}{\partial h}\Big|_{(\beta,\lambda,0)}.
%\end{align*}
Therefore, this derivative is equal to $0$ and since $\Phi_N$ is convex in $h$, $\Phi_N$
is non-decreasing in $h$. Then, the convergence of $\Phi_N$ to $\Phi$ implies
that $\Phi$ is also non-decreasing in
$h$.
%We let $n$ and $v$ tend to $\infty$ in \eqref{V_{1}} and we add $\lambda h$ on both sides.
The step $2$ gives,
for $h>0$, that
%\begin{equation}\label{cash}
$\limsup_{a \to 0} \frac{1}{a^{2}} \Phi(a \beta, a \lambda, a h)\leq
\tilde{\Phi}(\beta \Sigma, \lambda, h)$.
%\end{equation}
Since $\Phi$ is non-decreasing in $h$, the
former inequality implies,
$\limsup_{a \to 0} \frac{1}{a^{2}} \Phi\left(a \beta, a \lambda, 0\right)\leq
\tilde{\Phi}(\beta \Sigma, \lambda, h)$. Then we let $h \downarrow 0$ and the proof of the convergence
\eqref{theo31}
is completed.

\subsubsection{Part 2: proof of the convergence \eqref{corol}}
In this section, we assume that Theorem \ref{theo} is satisfied.
%We recall that in the case $\beta=0$
%the corollary is proven in \cite{BDH}, therefore we take $\beta\neq 0$.
%First we consider the case $\beta \Sigma<\beta_c$, which means that $x\to K_c^x$
%is continuous around $K_c^{\beta\Sigma}$.
We consider $\beta\neq 0$. We prove the convergence \eqref{corol} by applying Theorem \ref{theo} with particular parameters. However we have to
take into account the fact that there may exist $j\in \{-K,\dots,K\}$
such that $\mathbb{E}(\gamma_{1}^{j})=0$. Therefore, as we did in Step $2$ of the proof of
\eqref{theo31}
we consider $\mu>0$ small enough, such that,
$\mathbb{E}(\gamma_{i}^{j}+\mu)\neq 0$ and $\mathbb{E}(\gamma_{i}^{j}-\mu)\neq 0$ for every $j\in \{-K,\dots,K\}$.
Then,
we use the result of Theorem \ref{theo} with the variables $\gamma_{i}^{j}+\mu$ for $i\geq 1$
and $j\in\{-K,\dots,K\}$. We denote by
$\Psi_{\mu}$ the associated excess free energy and we let
$\Sigma_{\mu}=\Sigma+(2K+1)\mu$. Then we denote
$\rho=1/n$, $\mu_{1}=\beta \Sigma_\mu+1/n$, $h=(1+2/n)K_{c}^{\mu_{1}}$, $h^{'}=K_{c}^{\mu_{1}}$,
$\beta_{1}=\beta_{2}=\beta$, and $\lambda=1$.
For $a$ small enough, the first inequality of Theorem \ref{theo}
gives
\begin{equation}\label{Arrr}
\textstyle{\frac{1}{a^{2}}}\Psi_\mu\Big(a\beta, a, a \big(1+\frac{2}{n}\big)K_{c}^{\beta \Sigma_\mu
+\frac{1}{n}}\Big)\leq
\big(1+\frac{1}{n}\big) \tilde{\Psi}
\Big(\beta \Sigma_\mu+\frac{1}{n}, 1, K_{c}^{\beta \Sigma_\mu+\frac{1}{n}}\Big).
\end{equation}
By definition of $K_{c}^{(.)}$, the right hand side of \eqref{Arrr} is equal to zero. Moreover
$\Psi\leq \Psi_\mu$ for $\mu>0$. Therefore,
we have the inequality
%\begin{equation}\label{arrr}
$\textstyle{\limsup_{a\to \infty}}\,h_{c}^{a\beta}(a)/a\leq (1+2/n)
K_{c}^{\beta\Sigma_\mu+1/n}$.
%\end{equation}
Then, we let $n\uparrow\infty$ and $\mu\downarrow 0$ and since $x\to K_{c}^{x}$ is continuous in $\beta\Sigma$,
the former inequality
becomes $\limsup_{a\to \infty}h_{c}^{a\beta}(a)/a \leq K_{c}^{\beta\Sigma}$.
It remains to prove the opposite inequality. To that aim, we apply the second inequality of Theorem \ref{theo}
with the variables $\gamma_i^j-\mu$ for $i\geq 1$ and $j\in\{-K,\dots,K\}$ and
with the parameters $\rho=1/n$, $\mu_{2}=\beta \Sigma_{-\mu}-1/n$, $h=K_{c}^{\mu_{2}}-1/n$,
$h^{'}=(K_{c}^{\mu_{2}}-2/n)/(1+1/n)$,
$\beta_{1}=\beta_{2}=\beta$, and $\lambda=1$. For $a$ small enough we obtain
\begin{equation}\label{Arrr1}
\textstyle{\tilde{\Psi}}\Big(\beta\Sigma_{-\mu}-\frac{1}{n}, 1, K_{c}^{\beta \Sigma_{-\mu}
-\frac{1}{n}}-\frac{1}{n}\Big)
\leq \frac{1+1/n}{a^{2}}\ \Psi_{-\mu}
\Big(a\beta, a, \frac{a}{1+1/n}\,\big(K_{c}^{\beta \Sigma_{-\mu}-\frac{1}{n}}-\frac{2}{n}\big)\Big).
\end{equation}
Therefore, since the l.h.s. of \eqref{Arrr1} is strictly positive and
since $\psi\geq \psi_{-\mu}$ for all $\mu>0$, we can write the inequality
$\liminf_{a\to \infty}h_{c}^{a\beta}(a)/a\geq (K_{c}^{\beta \Sigma_{-\mu}-1/n}-2/n)/(1+1/n)$.
%\end{equation}
Finally, by continuity of $x\mapsto K_c^{x}$ around $\beta\Sigma$, we let $n\uparrow\infty$ and $\mu\downarrow 0$
and it completes the proof
of \eqref{corol}.% in the case $\beta \Sigma<\beta_c$.

As in the step $2$ of the proof of \eqref{theo31}, the case $\gamma_i^j=0$
for all $i\geq 1$ and all $j\in\{-K,\dots,K\}$
gives us directly \eqref{corol} in the case $\beta=0$.

%Now, we consider the case $\beta \Sigma>\beta_c$, and we let $M$ be a constant $>0$.
%We consider again the variables $\gamma_i^j-\mu$ for $i\geq 1$ and $j\in\{-K,\dots,K\}$ but we take $\mu$
%small enough such that $\beta \Sigma_{-\mu}=\beta \Sigma-(2K+1)\mu>\beta_c$. Then, with
%the parameters $\rho=1/n$, $\mu_{2}=\beta \Sigma_{-\mu}-1/n$, $h=M(1+2/n)$, $h^{'}=M$,
%$\beta_{1}=\beta_{2}=\beta$, and $\lambda=1$ the second inequality of Theorem \ref{theo} gives for $a$ small enough
%\begin{equation}\label{Arrr2}
%\textstyle{\tilde{\Psi}}\Big(\beta \Sigma_{-\mu}-\frac{1}{n}, 1, M (1+\frac{2}{n})\Big)
%\leq \frac{1+1/n}{a^{2}}\ \Psi_{-\mu}
%\Big(a\beta, a,a M\Big).
%\end{equation}
%Since $\beta \Sigma_{-\mu}>\beta_c$ we can choose $n$ large enough such that
%$\beta \Sigma_{-\mu}-\frac{1}{n}>\beta_c$. This implies that the
%l.h.s. of \eqref{Arrr2} is strictly positive and since $\Psi\geq \Psi_{-\mu}$ we obtain
%$\liminf_{a\to \infty}h_{c}^{a\beta}(a)/a\geq M$. This is true for every $M>0$, therefore
%$\lim_{a\to \infty}h_{c}^{a\beta}(a)/a=\infty$ and the proof of the corollary is completed.

\subsection{Proof of Theorem \ref{theo}}\label{sec33}

%As defined above, in all this proof $\left(S_{N}\right)_{N\geq 0}$ is a simple random walk.
%We first enounce a definition which will be very useful in the proof of the theorem.

%\subsubsection{Coarse graining}

%We can now perform the proof of theorem 3, but for simplicity
% instead of looking at disorder spread on the layer of width $2K$ around the origin, we
% only consider a disorder spread
%on the $x$ axis and of expectation $1$ (namely $K=0$, $E(w_{1}^{0})=1$, then $\Sigma=1$). In that way we avoid
%complicated expressions of the Hamiltonian
%without switching any difficulties because the steps 1,3,4 and 5 would be the same with the $2K$ large disorder.
% The difference would appear at step 2, when we apply the lemma, but we have proven it in complete generality
%at section $5$
%(with the sum of
%$2K$ local times), so step 2 works also with the multi-interface disorder.
\begin{rem}
\rm{We will only consider in this proof the case $\beta_1>0$ and $\beta_2>0$. Indeed, if for instance
$\beta_1<0$, we transform all the variables $(\gamma_i^j)_{\{i\geq 1,j\in I_1\}}$ into
$(-\gamma_i^j)_{\{i\geq 1,j\in I_1\}}$ and we take $-\beta_1$ instead of $\beta_1$.}
\end{rem}

First, we define a relation (previously introduced in \cite{BDH}), which is very useful to carry out the proof.
\begin{definition}\label{def11}
let $f_{t,\gep,\delta}(a,h,\beta_{1},\beta_{2})$ and $g_{t,\gep,\delta}(a,h,\beta_{1},\beta_{2})$
be real-valued functions.
The relation $f<<g$ occurs if for every $\beta_{3}>\beta_{1}$, $\beta_{2}>\beta_{4}$, $\rho>0$, and $h>h^{'}\geq 0$
satisfying $(1+\rho)h^{'}<h$,
there exists $\delta_{0}$ such that for $0<\delta<\delta_{0}$ there exists $\gep_{0}(\delta)$ such that
for $0<\gep<\gep_{0}$ there exists $a_{0}(\gep, \delta)$ satisfying %$ a<a_{o}(\gep, \delta)$
\begin{multline}
\limsup_{t \to \infty} f_{t,\gep,\delta}(a,h,\beta_{1},\beta_{2})-
(1+\rho) g_{t(1+\rho)^{2},\gep(1+\rho)^{2},\delta(1+\rho)^{2}}(a(1+\rho),h^{'},\beta_{3},\beta_{4})\leq 0\\
for \ \ 0<a<a_{0}
\end{multline}
\end{definition}
In this proof we consider some functions of the form
\begin{equation*}
F_{t,\gep,\delta}(a,h,\beta_{1},\beta_{2})=\mathbb{E}\Big[\frac{1}{t}\log
E\big(\exp(a H_{t,\gep,\delta}(a,h,\beta_{1},\beta_{2}))\big)\Big]
,\end{equation*}
and we denote
\begin{itemize}
\item $F^{1}_{t,\gep,\delta}(a,h,\beta_{1},\beta_{2})=\frac{1}{a^{2}} \Psi_{\lfloor t/a^{2}\rfloor}(a\beta_{1}, a\beta_{2},a,a h)$\\
\item ${F}^{7}_{t,\gep,\delta}(a,h,\beta_{1},\beta_{2})= \tilde{\Psi}_{t}(\beta_{1}\Sigma_{1}
+\beta_{2}\Sigma_{2},1, h).$
\end{itemize}
The proof of \eqref{eqgen} will consist in showing that $F^{1}<<F^{7}$ and $F^{7}<<F^{1}$ (denoted by
$F^{1}\sim F^{7}$). To that
aim, we will create the intermediate functions $F_{2},\dots,F_{6}$ associated with slight modifications
of the Hamiltonian to transform, step by step, the discrete Hamiltonian into the continuous one. As the relation
$\sim$ is transitive, we will prove at every step that $F^{i}\sim F^{i+1}$, to conclude finally
that $F^{1}\sim F^{7}$.

\subsection{Scheme of the proof}\label{shema}

To show that $F^{i}<<F^{i+1}$ we let $H^{i}=H^{I}+H^{II}$ and, by the H\"older inequality, we can bound  $F^{i}$
from above as follows
\begin{multline*}
F^{i}_{t,\gep,\delta}(a,h,\beta)\leq \textstyle{\frac{1}{t(1+\rho)}} \mathbb{E}\big[\log E\big(\exp(a (1+\rho)
H^{I})\big)\big]
+ \textstyle{\frac{1}{t(1+\rho^{-1})}} \mathbb{E}\big[\log E\big(\exp(a (1+\rho^{-1}) H^{II})\big)\big].
\end{multline*}
Thus, if we choose $H^{I}=H^{i+1}_{t(1+\rho)^{2},\gep(1+\rho)^{2},\delta(1+\rho)^{2}}(a(1+\rho),h^{'},\beta_{3}
,\beta_{4}),$
we obtain
\begin{multline*}
F^{i}_{t,\gep,\delta}(a,h,\beta_{1}
,\beta_{2}))-(1+\rho)F^{i+1}_{t(1+\rho)^{2},\gep(1+\rho)^{2},\delta(1+\rho)^{2}}
(a(1+\rho),h^{'},\beta_{3},\beta_{4}))\\ \leq \textstyle{\frac{1}{t(1+\rho^{-1})}}
\mathbb{E}\big[\log E\big(\exp(a (1+\rho^{-1}) H^{II})\big)\big]
.\end{multline*}
Then, it suffices to prove that
$\limsup_{t \to \infty} 1/t \log \mathbb{E}E\left(\exp\left(a (1+\rho^{-1}) H^{II})\right)\right)\leq 0$ for
$a,\epsilon$ and $\delta$ small enough.

We can assume without problem that
$\gep/a^2$, $\delta/\gep$ and $t/\delta$ are all
integers. In this way we avoid the brackets in the formulas.

\subsection{Step 1}

The first Hamiltonian that we consider in this proof is given by
\begin{equation*}
H_{t,\gep,\gd}^{(1)}(a,h,\gb_{1},\gb_{2})
=-2\sum_{i=1}^{t/a^2}\Delta_{i} (w_{i}+ah)+\,\gb_{1} \sum_{j\in I_{1}  }\sum_{i=1}^{t/a^2}\gamma_{i}^{j}
\ind_{\{S_{i}=j\}}+\,\gb_{2} \sum_{j\in I_{2}}\sum_{i=1}^{t/a^2}\gamma_{i}^{j}
\,\ind_{\{S_{i}=j\}}
,\end{equation*}
 with $\Delta_{i}=1$ if $\Lambda_{i}=-1$ and  $\Delta_{i}=0$ if $\Lambda_{i}=1$.\\
We define some notation to build the intermediate Hamiltonians (see Fig. $3$).
\begin{itemize}
\item $\gs_{0}=0,\,\ \ i^{v}_{0}=0$\ \ and\ \
$i^{v}_{k+1} = \inf \,\{\,n>\gs_{k}\varepsilon/a^2+\gd/a^2 : S_{n}=0 \}$
\item $m\ =\ \inf\,\{k\geq1 : i^v_{k}>t/a^2\}$
\item  $i_{k}=i_{k}^{v}$\ \ for\ \ $k<m$ and\ \ $i_{m}=t/a^{2}$
\item $\gs_{k+1}=\ \inf\,\{\,n\geq 0 : i_{k+1}\in\ ](n-1)\varepsilon/a^2,n \varepsilon/a^2]\}$,
\item  \ $\overline{I}_{k}\ =\
](\gs_{k-1})\,\varepsilon/a^2,\gs_{k}\,\varepsilon/a^2\,] \
\cap\ ]0,t/a^2]$,\ \
$s_{k+1}=\ \sign\,\Delta_{i_{k+1}-1}.$
\end{itemize}
%figure explicative
%We give an example of this construction with the Fig.$3$.

\vspace{2cm}
\setlength{\unitlength}{0.35cm}
\begin{picture}(15,5)(0,0)

\put(26.8,4.4){\vector(-1,-4){1}}
\put(26,5.5){$\sigma_{m-1}\frac{\epsilon}{a^{2}}$}

\put(23.5,-6){\vector(1,4){1.5}}
\put(21.6,-6.7){$i_{m-1}$}

\put(-1,4.5){\vector(1,-4){1}}
\put(10.2,-4){\vector(1,4){1}}
\put(9.2,-5.3){$i_{1}$}
%\put(6.3,-4.5){\vector(1,4){1}}
%\put(4.2,4.5){\vector(-1,-4){1}}
\put(9.4,4.3){\vector(1,-4){1}}
\put(13.5,4.3){\vector(-1,-4){1}}
\put(38,5.5){$\sigma_{m}\frac{\epsilon}{a^{2}}$}
\put(39.3,4.4){\vector(-1,-4){1}}
\put(38,-4){\vector(-1,4){1}}
\put(37.6,-5.3){$i_{m}$}
%\put(27,-4.5){\vector(-1,4){1}}
% {\thicklines
%  \qbezier(0,5)(3,5)(6,5)
% }
% \put(4,5.5){$ i_{k-1}$}

%\put(-1,4.8){$1$}
\put(25.7,0){\circle*{.35}}
% \put(17,-6){$(\sigma_{m-1}-1)\frac{\epsilon}{a^{2}}$}
\put(24.2,0){\circle*{.35}}
%\put(20.9,-4.3){\vector(3,4){3}}
\put(36.5,0){\circle*{.35}}
\put(38.1,0){\circle*{.35}}
{\qbezier(35,0.4)(35,0)(35,-0.4)}
\put(34.5,2){$\frac{t}{a^{2}}$}
{\qbezier(0,0)(0.2,0.2)(0.6,0.6)}
{\qbezier(0.6,0.6)(1.2,0)(1.6,-0.4)}
{\qbezier(1.6,-0.4)(2.2,0.2)(2.9,0.9)}
{\qbezier(2.9,0.9)(3.1,0.7)(3.3,0.5)}
{\qbezier(3.3,0.5)(4,1.2)(4,1.2)}
{\qbezier(4,1.2)(6,-0.8)(6,-0.8)}
{\qbezier(6,-0.8)(7.5,0.7)(7.5,0.7)}
{\qbezier(7.5,0.7)(8,0.2)(8,0.2)}
{\qbezier(8,0.2)(9,1.2)(9.5,1.7)}
{\qbezier(9.5,1.7)(11.5,-0.3)(11.5,-0.3)}
{\qbezier(11.5,-0.3)(12.8,1)(12.8,1)}
{\qbezier(12.8,1)(14.2,-1)(14.2,-1)}
{\qbezier(14.2,-1)(14.7,-0.5)(14.7,-0.5)}
{\qbezier(22.5,-0.5)(23,-1)(23.5,-1.5)}
{\qbezier(23.5,-1.5)(25.5,0.5)(26.5,1.5)}
{\qbezier(26.5,1.5)(27,1)(27.5,0.5)}
{\qbezier(27.5,0.5)(28,1)(29,2)}
{\qbezier(29,2)(29.5,2.5)(30,3)}
{\qbezier(30,3)(33,0)(34,-1)}
{\qbezier(34,-1)(34.3,-0.7)(34.5,-0.5)}
{\qbezier(34.5,-0.5)(34.8,-0.8)(35.5,-1.5)}
{\qbezier(35.5,-1.5)(36.5,-0.5)(37.5,0.5)}
{\qbezier(37.5,0.5)(38,0)(38.5,-0.5)}

{\qbezier(0,0)(5,0)(15,0)}
{\thicklines
\qbezier[14](15,0)(18,0)(21,0)
}
\put(21,0){\vector(1,0){19}}
\put(7.5,0){\circle*{.35}}
\put(0,0){\circle*{.35}}
%{\qbezier(0.1,-0.2)(3.2,-4)(7.4,-0.2)}
\put(10.6,0){\circle*{.35}}
\put(12.4,0){\circle*{.35}}
%{\qbezier(25.8,-0.2)(29,-4)(33.4,-0.2)}
\put(33.5,0){\circle*{.35}}
\put(29,-4.5){$\frac{\delta}{a^{2}}$}
\put(3,-4.5){$\frac{\delta}{a^{2}}$}
%\put(4,-6){$\sigma_{0}\frac{\epsilon}{a^{2}}+\frac{\delta}{a^{2}}$}
\put(-3.5,6){$\sigma_{0} \frac{\gep}{a^{2}}$}
\put(5.5,5.5){$(\sigma_{1}-1) \frac{\gep}{a^{2}}$}
\put(12.8,5.5){$\sigma_{1} \frac{\gep}{a^{2}}$}
\put(0,8){\small
             Fig.\ 3:}
             {\thicklines
             {\qbezier[6](0,0)(0,-1.5)(0,-3)}}
             {\thicklines
             {\qbezier[6](7.5,0)(7.5,-1.5)(7.5,-3)}}
\put(0.2,-3){\vector(1,0){7.2}}
\put(7.3,-3){\vector(-1,0){7.2}}
{\thicklines
             {\qbezier[6](25.7,0)(25.7,-1.5)(25.7,-3)}}
             {\thicklines
             {\qbezier[6](33.5,0)(33.5,-1.5)(33.5,-3)}}
\put(26.2,-3){\vector(1,0){7.2}}
\put(33,-3){\vector(-1,0){7.2}}
\end{picture}

\vspace{3.5cm}
We define the first transformation of the Hamiltonian
\begin{align*}
H_{t,\gep,\gd}^{(2)}(a,h,\gb_{1},\gb_{2})=-2\sum_{k=1}^{m}s_{k}\bigg[\sum_{i \in \overline{I}_{k}}w_{i}+a h|
\overline{I}_{k}|\bigg]+\sum_{i=1}^{t/a^2}
\gb_{1}\sum_{j\in I_{1} }\gamma_{i}^{j}
\ind_{\{S_{i}=j\}}
 +\gb_{2}\sum_{j\in I_{2}}\gamma_{i}^{j}\ind_{\{S_{i}=j\}}
\end{align*}
and we want to show that $F_{1}<<F_{2}$. To that aim, we denote
\begin{align}\label{pat}
\nonumber \textstyle{H}^{II}=\,-&\textstyle{2}\sum_{i=1}^{t/a^2}\,\Delta_{i} (w_{i}+ah)\,
+2\sum_{k=1}^{m}s_{k}\Big(\sum_{i \in \overline{I}_{k}}
w_{i}+a (1+\rho)h^{'}|
\overline{I}_{k}|\Big)\, \\+&
\textstyle{(\beta_{1}-\beta_{3})}\sum_{j\in I_{1} }\sum_{i=1}^{t/a^2}\gamma_{i}^{j}
\,\ind_{\{S_{i}=j\}}+(\beta_{2}-\beta_{4})\sum_{j\in I_{2}}\sum_{i=1}^{t/a^2}\gamma_{i}^{j}
\,\ind_{\{S_{i}=j\}}
,\end{align}
and it remains to prove that
%\begin{equation}\label{pat3}
$\limsup_{t\rightarrow \infty}\frac{1}{t} \log
E\mathbb{E}(\exp(a(1+\rho^{-1})H^{II}))\leq 0$.
%\end{equation}
We integrate over the disorder $\gamma$ and
the third and forth terms of the right hand side of \eqref{pat} give some contributions of the form
\begin{align}\label{pat4}
\nonumber \textstyle{\exp}\big(\sum_{j\in I_{p}}\,\sum_{i=1}^{t/a^{2}}\log \mathbb{E}
&\big[\exp\big((\beta_{p}-\beta_{2+p})a(1+\rho^{-1})
\gamma_{i}^{j}\big)\big]
\,\ind_{\{S_{i}=j\}}\big) \quad\text{for $p=1$ and $p=2$}.
%\times \exp\bigg(\sum_{j\in I_{2}}\sum_{i=1}^{t/a^{2}}\log \mathbb{E}
%\big((\beta_{2}-\beta_{4})a(1+\rho^{-1})
%\gamma_{i}^{j}\big)
%\,\ind_{\{S_{i}=j\}}\bigg).
\end{align}
Since $\mathbb{E}(\exp(\lambda |\gamma_{1}^{j}|))<\infty$ for every $j\in\{-K,\dots,K\}$ and
$\lambda>0$, we can write a
first order Taylor expansion of $\log\mathbb{E}(\exp(A a \gamma_{1}^{j}))$ when $a \downarrow 0$. It gives
\begin{equation}\label{Fd}
\log\mathbb{E}\big(\exp\big(A a \gamma_{1}^{j}\big)\big)=A a \mathbb{E}\big(\gamma_{i}^{j}\big)+o(a).
\end{equation}
We assume in this proof that $\mathbb{E}(\gamma_{1}^{j})\neq 0$ for every
$j\in\{-K,\dots,K\}$ (see the assumptions of Theorem \ref{theo}) and therefore $\{-K,\dots,K\}=I_{1}\cup I_{2}$.
For every $i \in I_{1}$, $\mathbb{E}(\gamma_{1}^{j})>0$, and $\beta_{1}-\beta_{3}<0$. Thus,
by \eqref{Fd}, we obtain, for $a$ small enough, that
\begin{equation}\label{strik}
\textstyle{\sum_{j\in I_{1} }}\sum_{i=1}^{t/a^{2}}\log \mathbb{E}
\big((\beta_{1}-\beta_{3})a(1+\rho^{-1})
\gamma_{i}^{j}\big)
\,\ind_{\{S_{i}=j\}}\leq 0.
\end{equation}
The sum over $I_{2}$ satisfies the same inequality for $a$ small enough because $\beta_{2}-\beta_{4}>0$ and
$\mathbb{E}\big(\gamma_{i}^{j})<0$ when $j \in I_{2}$.
Therefore, we can remove the third and forth terms of $H^{II}$ in \eqref{pat} and
by rewriting $\sum_{i=1}^{t/a^{2}}$ as $\sum_{k=1}^{m_{t/a^2}}\sum_{i\in \overline{I}_{k}}$,
we can rewrite $H^{II}$ as
\begin{align*}
H^{II}\hspace{-0.1cm}=\hspace{-0.1cm}-2\hspace{-0.2cm}
\sum_{k=1}^{m_{t/a^2}}\hspace{-0.1cm}\sum_{i\in \overline{I}_{k}}
w_{i}\left(\Delta_{i}-s_{k}\right)\hspace{-0.05cm}-\hspace{-0.05cm}
2a (1+\rho) h^{'}\hspace{-0.2cm}\sum_{k=1}^{m_{t/a^2}}
\hspace{-0.1cm}\sum_{i\in \overline{I}_{k}}
\left(\Delta_{i}-
s_{k}\right)
\hspace{-0.05cm}-\hspace{-0.05cm}2a (h-(1+\rho) h^{'})\hspace{-0.2cm}\sum_{k=1}^{m_{t/a^2}}\hspace{-0.1cm}\sum_{i\in \overline{I}_{k}}
\Delta_{i}
%\sum_{i\in
%\overline{I}_{k}}
%\Delta_{i}+\left(\beta-\beta'\right)\sum_{i=1}^{t/a^2}\,w_{i}\ind{\{S_{i}=0\}}.
\end{align*}
Thus, we integrate over the disorder $w$ which is independent of the random walk. But, since $\mathbb{E}(w_{i})=0$
and $\mathbb{E}(\exp(\lambda |w_{1}|))<\infty$ for every
$\lambda >0$, a second order expansion gives that
for every $c\in \mathbb{R}$ there exists $A>0$ such that for $a$ small enough
\begin{equation}
\log\mathbb{E}\left(\exp(c\, a\, w_{i}\, \left(\Delta_{i}-s_{k}\right)\right))\leq A a^{2} |\Delta_{i}-s_{k}|
.\end{equation}
Finally, we have to prove, for $A>0$ and $B>0$ and for $\delta,\gep,a$ small in the sense of Definition \ref{def11}, that
\begin{equation}\label{eq7}
\textstyle{\limsup_{t\to \infty}}\frac{1}{t} \log E\Big[\exp\big(A a^{2}\sum_{k=1}^{m_{t/a^2}}\sum_{i\in
\overline{I}_{k}}\mid s_{k}-\Delta_{i}\mid -B a^{2}\sum_{i=1}^{t/a^{2}}\Delta_{i}\big)\Big]\leq 0
.\end{equation}
This is explicitly proven in \cite{BDH} (page $1355$), and completes the Step 1 because
the proof of $F_{2}<<F_{1}$ is very similar and consists essentially in showing \eqref{eq7}.

\subsection{Step 2}

In this step we aim at transforming the disorder $w$ into a sequence $(\hat{w}_{i})_{i\geq1}$ of independent
random variables of law $\mathcal{N}_{0,1}$. To that aim, we use a coupling method developed in \cite{QMS}
to define on the same probability space and for every $j \in \mathbb{N}\setminus\{0\}$ the variables $(w_{i})_{i \in \{(j-1) \gep
/a^{2}+1,\dots,j \gep/a^{2}\}}$ and
some independent variables of law $\mathcal{N}_{0,1}$, denoted by $(\hat{w}_{i})_{i \in \{(j-1) \gep
/a^{2}+1,\dots,j \gep/a^{2}\}}$, such that for every $p>2$ and $x>0$
\begin{equation}\label{KMT}
\mathbb{P}\bigg( \Big| \sum_{i=(j-1)\gep/a^{2}+1}^{j\gep/a^{2}}w_{i}-\hat{w}_{i}\Big| \geq x\bigg) \leq
\textstyle{\frac{(Ap)^{p}\  \gep}{x^{p}\  a^{2}}\ \mathbb{E}\left(w_{1}^{p}\right)}
.\end{equation}
These constructions are made independently on all blocs $\{(j-1) \gep/a^{2}+1,\dots,j \gep/a^{2}\}$.
Thus, we can form the third Hamiltonian as follow
\begin{align*}
H_{t,\gep,\gd}^{(3)}(a,h,\gb_{1},\gb_{2})\, =-2\sum_{k=1}^{m}s_{k}
\Bigg[\sum_{i \in \overline{I}_{k}}\hat{w}_{i}+a h|
\overline{I}_{k}|\Bigg]+
\sum_{i=1}^{t/a^2}
\gb_{1}\sum_{j\in I_{1}}\gamma_{i}^{j}\ind_{\{S_{i}=j\}}
+\gb_{2}\sum_{j\in I_{2}}\gamma_{i}^{j}
\ind_{\{S_{i}=j\}}.
\end{align*}
To prove that $F^{2}<<F^{3}$, we need the Hamiltonian $H^{II}$. It takes the value
\begin{equation}
H^{II}=H_{t,\gep,\gd}^{(2)}(a,h,\gb_{1},\gb_{2})-H_{t(1+\rho)^{2},\gep(1+\rho)^{2},\gd(1+\rho)^{2}}^{(3)}
(a(1+\rho),h^{'},\gb_{3},\gb_{4})
.\end{equation}
As in Step 1 (see \eqref{strik}) we delete the two pinning terms in $H^{II}$
%$$\left(\beta_{1}-\beta_{3}\right)\sum_{j\in I_{1}}\sum_{i=1}^{t/a^2}
%\gamma_{i}^{j}\,\ind_{\{S_{i}=j\}}+\left(\beta_{2}-\beta_{4}\right)\sum_{j\in I_{2}}\sum_{i=1}^{t/a^2}
%\gamma_{i}^{j}\,\ind_{\{S_{i}=j\}},$$
and it is sufficient  to consider
\begin{align*}
H^{II}&=-2 \sum_{k=1}^{m}s_{k}\sum_{i\in \overline{I}_{k}}(w_{i}-\hat{w}_{i})+2a \sum_{k=1}^{m}s_{k}(h-(1+\rho)h')|
\overline{I}_{k}|\\
&\leq 2\sum_{k=1}^{m}s_{k}\bigg(\sum_{j=\sigma_{k-1}+1}^{\sigma_{k}} \bigg|\sum_{i=j\gep/a^{2}+1}^
{(j+1)\gep/a^{2}}w_{i}-\hat{w}_{i}\bigg|-(h-(1+\rho)h')\textstyle{\frac{\gep}{a}}\bigg)
.\end{align*}
We want to prove that
$\limsup_{t \to \infty} 1/t \log \mathbb{E}E\left(\exp\left(a (1+\rho^{-1}) H^{II})\right)\right)\leq 0$.
By independence of $(w,\hat{w})$ on each blocs $\{(j-1) \gep/a^{2}+1,\dots,j \gep/a^{2}\}$,
it suffices
to show that for every $C>0$ and $B>0$
\begin{equation}\label{main}
\textstyle{\mathbb{E}}\Big[\exp\big(C a  \big|\sum_{i=1}^{\gep/a^{2}}
w_{i}-\hat{w}_{i}\big|-B \gep\big)\Big]\leq 1\  \text{for\  $\gep$,\ \text{and}\ $a$\  small enough.}
\end{equation}
We prove this point as follows,
\begin{align}\label{eqi}
\textstyle{\mathbb{E}}\Big[\exp\big(C a  \big|\sum_{i=1}^{\gep/a^{2}}
w_{i}-\hat{w}_{i}\big|\big)\Big]\leq \sum_{k=N}^{+\infty}e^{Ca (k+1) \frac{\gep}{\sqrt{a}}}
\ \,\mathbb{P}\Big(\big|\sum_{i=1}^{\gep/a^{2}}
w_{i}-\hat{w}_{i}\big|\geq k \textstyle{\frac{\gep}{\sqrt{a}}}\Big)+ e^{C N \sqrt{a} \gep}
.\end{align}
By using \eqref{KMT} and the fact that $\mathbb{E}(w_{1}^{k})\leq R^{k}$, we obtain that for every $j$ and $k$
$\geq 1$
\begin{equation}\label{KMA}
\mathbb{P}\bigg(\Big| \sum_{i=(j-1)\gep/a^{2}+1}^{j\gep/a^{2}}w_{i}-\hat{w}_{i}\Big| \geq \textstyle{\frac{k \gep}
{\sqrt{a}}}\bigg) \leq
\frac{(AR\sqrt{a})^{k}}{\gep^{k-1}\, a^{2}}
.\end{equation}
We consider \eqref{eqi} with $N=5$, and we use \eqref{KMA} to obtain
\begin{align}\label{eqo}
\nonumber \textstyle{\mathbb{E}\Big[\exp\big(C a  \big|\sum_{i=1}^{\gep/a^{2}}
w_{i}-\hat{w}_{i}\big|\big)\Big]\leq e^{5 C \sqrt{a} \gep}
+\gep\ \frac{e^{C\sqrt{a}\gep}}{a^{2}}\, \sum_{k=5}^{+\infty}\,
\big(e^{C \sqrt{a} \gep}\textstyle{\frac{AR\sqrt{a}}{\gep}}\big)^{k}}.
%\leq & \exp(5 A \sqrt{a} \gep)+\frac{(A R)^{5} \sqrt{a}\exp
%\left(6 A\sqrt{a}\gep\right)}{\gep^{4} \left(1-\exp(A \sqrt{a} \gep) A R \sqrt{\gep}\right)}
\end{align}
Therefore, for $\gep>0$ fixed, there exists $K(\gep,a)>0$ which tends to zero when $a$ tends to zero, and satisfies
\begin{equation*}
\textstyle{\mathbb{E}}\Big[\exp\big(C a  \big|\sum_{i=1}^{\gep/a^{2}}
w_{i}-\hat{w}_{i}\big|\big)\Big]\leq (1+K(\gep,a))\ \ e^{5 C \gep \sqrt{a}}
.\end{equation*}
This implies \eqref{main}, and completes the Step $2$ because the proof of $F^{3}<<F^{2}$ is exactly the same.

\subsection{Step 3}

In this step, we make a link between the discrete and the continuous models. For that, we take into
account the number of returns to the origin of the random walk, and the
local time of the Brownian motion.
We define, independently
of the random walk, an i.i.d. sequence $(l_{1}^{k})_{k\geq0}$
of local times spent in $0$ by a Brownian motion between $0$ and $1$. The law of this sequence is denoted by
$\chi$. Then, we build the new Hamiltonian
\begin{equation}\label{H3}
\textstyle{H}_{t,\gep,\gd}^{(4)}(a,h,\gb_{1},\gb_{2}) =\, -2\sum_{k=1}^{m}s_{k}\Big(\sum_{i \in \overline{I}_{k}}
\hat{w}_{i}+a h|\overline{I}_{k}|\Big)\, +
\frac{(\gb_{1}\Sigma_{1}+\gb_{2}\Sigma_{2})\sqrt{\delta}}{a}\sum_{k=1}^{m} l_{1}^{k}
.\end{equation}
As usual, to prove that $F_{3}<<F_{4}$, we consider $H^{II}$, in which we can already remove the term
$-2a(h-(1+\rho)h')\sum_{k=1}^{m} s_{k}\mid
\overline{I}_{k}\mid$ because it is negative. Therefore we can bound $H^{II}$ from above as follows
\begin{align*}
H^{II}%= -2a(h-&(1+\rho)h')\sum_{k=1}^{m} s_{k}\mid
%\overline{I}_{k}\mid+
%\beta_{1} \sum_{k=1}^{m}\ \sum_{j\in I_{1}}\sum_{\ i=i_{k-1}+1}^{i_{k}}
%\gamma_{i}^{j}\ind_{\{S_{i}=j\}}\\
%&+\beta_{2} \sum_{k=1}^{m}\ \sum_{j\in I_{2}}\sum_{\ i=i_{k-1}+1}^{i_{k}}
%\gamma_{i}^{j}\ind_{\{S_{i}=j\}}
%-\frac{(\beta_{4}\Sigma_{2}+\beta_{3}\Sigma_{1})\sqrt{\delta}\,\Sigma}{a}\sum_{k=1}^{m}l_{1}^{k}\\
\leq \beta_{1} \sum_{k=1}^{m}\ \sum_{j\in I_{1}} \sum_{\ i=i_{k-1}+1}^{i_{k}}&
\gamma_{i}^{j} \ind_{\{S_{i}=j\}}-\frac{\beta_{3}\Sigma_{1}\sqrt{\delta}}{a}\sum_{k=1}^{m}l_{1}^{k}\\
&+\beta_{2} \sum_{k=1}^{m}\ \sum_{j\in I_{2}} \sum_{\ i=i_{k-1}+1}^{i_{k}}
\gamma_{i}^{j} \ind_{\{S_{i}=j\}}-\frac{\beta_{4}\Sigma_{2}\sqrt{\delta}}{a}\sum_{k=1}^{m}l_{1}^{k}.
\end{align*}
To prove that
$\limsup_{t\rightarrow \infty} \frac{1}{t} \log
E_{P\otimes \chi}\mathbb{E}(\exp(a(1+\rho^{-1})H^{II}))\leq 0$,
we first apply the H\"older inequality (with the coefficients $p=q=2$),
and then we integrate over the disorder $\gamma$. Therefore, it remains to prove
for $x=1\ \text{and}\ 2$ that
%\begin{multline*}
%\limsup_{t \to \infty}\frac{1}{t}\log E_{P\otimes M}\mathbb{E}\bigg[\exp\bigg(
%2\sum_{k=1}^{m}\sum_{j=I_{x}}\sum_{\ i=i_{k-1}+1}^{i_{k}}\\
%a \beta_{x} \left(1+\rho^{-1}\right) \gamma_{i}^{j}
%\ind_{\{S_{i}=j\}}
%-2\beta_{x+2} \sqrt{\delta} \Sigma_{x} (1+\rho^{-1})l_{1}^{k}
%\bigg)\bigg]\leq 0.
%\end{multline*}
%We integrate over the disorder
%$\gamma$ and it remains to prove that for $x=1\ \text{and}\ 2$
\begin{align}\label{Kai}
\nonumber \limsup_{t \to \infty}\frac{1}{t}\log& E_{P\otimes \chi}\bigg[\exp\bigg(
\sum_{k=1}^{m}\sum_{j\in I_{x}}\sum_{\ i=i_{k-1}+1}^{i_{k}}\\
&\log \mathbb{E}\left(\exp\left(2 a \beta_{x} \left(1+\rho^{-1}\right) \gamma_{i}^{j}
\right)\right)
\ind_{\{S_{i}=j\}}
-2 \beta_{x+2} \sqrt{\delta} \Sigma_{x} (1+\rho^{-1})l_{1}^{k}
\bigg)\bigg]\leq 0.
\end{align}
For simplicity, in what follows we will use $E$ instead of $E_{P\otimes \chi}$.
We begin with the proof of \eqref{Kai} in the case $x=1$. To that aim,
we recall \eqref{Fd}, that gives
\begin{equation}
\log\mathbb{E}\left(\exp\left(2 a \beta_{1} \left(1+\rho^{-1}\right) \gamma_{i}^{j}
\right)\right)=2\mathbb{E}(\gamma_{1}^{j})a \beta_{1} \left(1+\rho^{-1}\right)+o(a).
\end{equation}
Therefore, we can choose $\beta''$ such that $\beta_{1}<\beta''<\beta_{3}$ and $a$ small enough to obtain for
every $j\in I_{1}$
the inequality
$\log\mathbb{E}(\exp(2 a \beta_{1} (1+\rho^{-1}) \gamma_{i}^{j}))\leq 2 a \beta''
(1+\rho^{-1}) \mathbb{E}(\gamma_{1}^{j})$.
Finally, since $\mathbb{E}(\gamma_{1}^{j})>0$ for every $j$, we can replace $(i_{k})_{k\in\{1,\dots,m\}}$
by $(i_{k}^{v})_{k\in\{1,\dots,m\}}$ (see the notation at the beginning of Step 1), and it remains to prove that
for $B>A>0$
\begin{equation}\label{in1}
\limsup_{t \to \infty}\frac{1}{t}\log E\bigg[\exp\bigg(
\sum_{k=1}^{m} \Big( A a\sum_{j\in I_{1}} \mathbb{E}(\gamma_{1}^{j})
\sum_{i=i^{v}_{k-1}+1}^{i_{k}^{v}}
\ind_{\{S_{i}=j\}}
-B \sqrt{\delta}\ \Sigma_{1}\  l_{1}^{k}\Big)\bigg)\bigg]\leq 0.
\end{equation}
For simplicity, we will use the notation $\mathbb{E}(\gamma_{1}^{j})=f(j)$, and
consequently $\Sigma_{1}=\sum_{j\in I_{1}}f(j)$.
%\sum_{i\in \{i_{k-1}+1,i_{k}\}}
For every $N$, we build a new filtration, i.e.,
$\textsl{F}_{N}=\sigma(A_{i^{v}_{N}}\cup \sigma(l_{0}^{1}
\dots,l_{1}^{N}))$ with $A_{k}=\sigma(X_{1},\dots,X_{k})$ and the
random variable
\begin{equation*}
M_{N}=\frac{\exp\big( \sum_{k=1}^{N}A a\sum_{j\in I_{1}}f(j)\ \sharp\{v\in \{i_{k-1}^{v}+1,i_{k}^{v}\}
: S_{v}=j\}-B\,\sqrt{\delta}\,\Sigma_{1}\,\sum_{k=1}^{N}l_{1}^{k}\big)}{\mu^{N} E\big(\exp
\big(A a \sum_{j\in I_{1}}f(j)\sharp\{
i \in \{0,\frac{\delta+\epsilon}{a^{2}}\} : S_{i}=j\}-B\sqrt{\delta}\,\Sigma_{1}\, l_{1}^{1}\big)\big)^{N}}
\end{equation*}
where $\mu$ is a constant $>1$. We will precise the value of $\mu$ later, to make sure that $M_{N}$ is a
positive super-martingale
with respect to $\left(F_{N}\right)_{N\geq 0}$.
%for every $k\in\mathbb{N}$ we have $\sharp\big\{j\in \{i_{k-1}+1,i_{k}\}
%: S_{j}=0\big\}\leq\sharp\big\{j\in \{i_{k-1},i_{k-1}+(\delta+\epsilon)/a^{2}\}
%: S_{j}=0\big\}$ (see Fig. $4$). We notice also that $m_{t/a^{2}}$ is a stopping time with respect
%to the filtration.
To that aim, for every $j\in \{-K,\dots,K\}$ we introduce
%\begin{equation*}
$P_{N}^{j}=\sharp\{u\in
\{i_{N-1}^{v}+1,i_{N}^{v}\}: S_{u}=j\}$,
%\end{equation*}
and we define the new filtration $(G_N)_{N\geq 1}$  by
%\begin{equation*}
$\textsl{G}_{N-1}=\sigma(\textsl{F}_{N-1}\cup \sigma
(X_{i_{N-1}^{v}+1},\dots,X_{i_{N-1}^{v}+
(\delta +\gep)/a^{2}},l_{1}^{N}))$.
%\end{equation*}
Then, we consider the quantity $E\left(M_{N}|\textsl{F}_{N-1}\right)$ and by independence of the
random walk excursions out of the origin we obtain
\begin{equation}\label{mart}
E(M_{N}|\textsl{F}_{N-1})=M_{N-1} \frac{\mu^{-1} E(\exp(A a\sum_{j\in I_{1}}f(j)\ P_{N}^{j}\,
-B\sqrt{\delta}\,\Sigma_{1}\,l_{1}^{N})\big|\textsl{F}_{N-1})}
{E(\exp(A a \sum_{j\in I_{1}} f(j)\
\sharp\{i \in \{0,\frac{\delta+\epsilon}{a^{2}}\} : S_{i}=j\}-B\sqrt{\delta}\Sigma_{1} l_{1}^{1}))}.
\end{equation}
We define $t_{N}=\inf\{i>i_{N-1}^{v}+(\delta+\gep)/a^{2} : S_{i}=0\}$ and notice that $t_{N}\geq i_{N}^{v}$
(see Fig. $4$ for an example in which $t_{N}>i_{N}^{v}$).

%%% BEGIN FIGURE 3: picture of $\kappa \mapsto \lambda_p(\kappa)$ %%%

\vspace{2cm}

\setlength{\unitlength}{0.35cm}

\begin{picture}(15,5)(0,0)

\put(0,0){\vector(1,0){40}}
\put(-1,-4.5){\vector(1,4){1}}
\put(9,-4.5){\vector(-1,4){1}}
\put(5.1,2.2){\vector(-1,-2){1}}
\put(38.2,2.2){\vector(-1,-2){1}}
\put(29.1,2.2){\vector(-1,-2){1}}
%\put(32,4.5){\vector(1,-4){1}}
 \put(27,-4.5){\vector(-1,4){1}}
% {\thicklines
%  \qbezier(0,5)(3,5)(6,5)
% }
\put(5,2.5){$ i_{N-1}$}
\put(29,2.5){$ i_{N}$}
\put(38,2.5){$ t_{N}$}
%\put(-1,4.8){$1$}
\put(24,-6){$\small\sigma_{N-1}\epsilon/a^{2}+\delta/a^{2}$}

{\qbezier(0,2)(0.5,1.5)(1,1)}
{\qbezier(1,1)(1.5,1.5)(2,2)}
{\qbezier(2,2)(4,0)(5,-1)}
{\qbezier(5,-1)(5.5,-0.5)(6,0)}
{\qbezier(6,0)(6.5,-0.5)(7.5,-1.5)}
{\qbezier(7.5,-1.5)(8,-2)(8.5,-2.5)}
{\qbezier(8.5,-2.5)(10.5,-0.5)(11.5,0.5)}
{\qbezier(11.5,0.5)(12,0)(12.5,-0.5)}
{\qbezier(12.5,-0.5)(13,0)(13.5,0.5)}
{\qbezier(13.5,0.5)(14.5,-0.5)(17,-3)}
{\qbezier(17,-3)(17.5,-2.5)(18,-2)}
{\qbezier(18,-2)(19,-3)(19.5,-3.5)}
{\qbezier(19.5,-3.5)(20.5,-2.5)(22.5,-0.5)}
{\qbezier(22.5,-0.5)(23,-1)(23.5,-1.5)}
{\qbezier(23.5,-1.5)(25.5,0.5)(26.5,1.5)}
{\qbezier(26.5,1.5)(27,1)(27.5,0.5)}
{\qbezier(27.5,0.5)(28,0)(29.5,-1.5)}
{\qbezier(29.5,-1.5)(31,0)(33,2)}
{\qbezier(33,2)(34,1)(34.5,0.5)}
{\qbezier(34.5,0.5)(35,1)(35.5,1.5)}
{\qbezier(34.5,0.5)(35,1)(35.5,1.5)}
{\qbezier(35.5,1.5)(37,0)(38.5,-1.5)}

{\thicklines
             {\qbezier[20](4,0)(4,3)(4,6)}}
             {\thicklines
             {\qbezier[20](32,0.1)(32,3)(32,6)}}
\put(6,6){\vector(1,0){25.9}}
\put(29,6){\vector(-1,0){24.9}}
\put(14.5,7.5){$\delta/a^{2}+\gep/a^{2}$}
\put(8,0){\circle*{.65}}
\put(0,0){\circle*{.65}}
\put(26,0){\circle*{.65}}
\put(8,-6){$\small\sigma_{N-1}\epsilon/a^{2}$}
\put(-4,-6){$\small\sigma_{N-1}\epsilon/a^{2}-\epsilon/a^{2}$}
\put(0,9.5){\small
             Fig.\ 4:}%:{\thicklines
   %\qbezier[6](3,0)(5,0)(7,0)}}

\end{picture}

\vspace{3cm}

%%%%%%% END FIGURE 3 %%%%%%%%%%%%%%%%%%%%%%%%%%%%%%%%%%%%%%
\noindent
Therefore,
we can write $P_{N}^{j}\leq B_{1,N}^{j}+B_{2,N}^{j}$ with
\begin{align}
B_{1,N}^{j}\hspace{-0.1cm}=\hspace{-0.15cm}\{v\hspace{-0.05cm}\in
\{i_{N-1}^{v}\hspace{-0.1cm}+\hspace{-0.1cm}1,\dots,i_{N-1}^{v}\hspace{-0.1cm}+
\hspace{-0.1cm}\textstyle{\frac{\delta+\gep}{a^{2}}}\}
\hspace{-0.05cm}:\hspace{-0.05cm} S_{v}\hspace{-0.1cm}=\hspace{-0.1cm}j\}\ \text{and}\
 B_{2,N}^{j}\hspace{-0.1cm}=\hspace{-0.15cm}\{v\hspace{-0.05cm}\in\hspace{-0.05cm}
\{i_{N-1}^{v}\hspace{-0.05cm}+\hspace{-0.05cm}\textstyle{\frac{\delta+\gep}{a^{2}}}
\hspace{-0.1cm}+1,\dots,t_{N}\}
\hspace{-0.05cm}:\hspace{-0.05cm} S_{v}\hspace{-0.1cm}=\hspace{-0.1cm}j\}.
\end{align}
We denote by $C$ the quantity $E[\exp(A a\sum_{j\in I_{1}}f(j)\ P_{N}^{j}\,
-B\sqrt{\delta}\,\Sigma_{1}\,l_{1}^{N})|F_{N-1}]$. Thus, since $B_{1,N}^{j}$ is measurable
with respect to $\textsl{G}_{N-1}$ and since $F_{N-1}\subset G_{N-1}$ we can write
\begin{align*}
%C&=E\Bigg[\exp\Big(A a\sum_{j\in I_{1}}f(j)\ P_{N}^{j}\,
%-B\sqrt{\delta}\,\Sigma_{1}\,l_{1}^{N}\Big)\Big|F_{N-1}\bigg]\\
C\leq E\bigg[\exp\Big(A a\sum_{j\in I_{1}}f(j)
\ B_{1,N}^{j}\,
-B\sqrt{\delta}\Sigma_{1}\,l_{1}^{N}\Big)\, E\bigg[\exp\Big(A a\sum_{j\in I_{1}}f(j)\ B_{2,N}^{j}\Big)
\Big|G_{N-1}\bigg]\ \bigg|F_{N-1}
\Bigg].
\end{align*}
We recall that $A_k=\sigma(X_{1},\dots,X_{k})$ and we let
$\Upsilon=E(\exp(A a\sum_{j\in I_{1}}f(j)\ B_{2,N}^{j})\,|G_{N-1})$. The fact that
the local times $(l_{1}^{1},\dots,l_{1}^{N})$ are independent of the random walk allows us to write the equality
$\Upsilon=E(\exp(A a\sum_{j\in I_{1}}f(j)\ B_{2,N}^{j})\,|A_{i^v_{N-1}+(\delta+\gep)/a^{2}
})$.
The strong Markov property can be applied here. In fact, if $(V_{n})_{n\geq 0}$ is a simple random
walk with $V_{0}=S_{i_{N-1}^{v}+(\delta+\gep)/a^{2}}$, and if $s=\inf\{n>1 : V_{n}=0\}$, we can write
\begin{equation*}
\textstyle{\Upsilon=E_{V}\big[\exp\big(A a\sum_{j\in I_{1}}f(j) \sharp\ \{i \in\{1,.,s\} :
V_{i}=j\}\big)\big]}
.\end{equation*}
Thus, if we denote $f=\max_{j\in I_{1}}\{f_{j}\}$, we can bound $\Upsilon$ from above as
\begin{equation}\label{HP}
\Upsilon\leq%E\Bigg(\exp\Bigg(A a\sum_{j=-K}^{K}f(j)\ B_{2}\Bigg)\,\bigg|G_{N-1}\Bigg)=
E_{V}\big[\exp\big(A a f \sharp\ \{i \in\{1,.,s\} :
V_{i}\in \{-K,\dots,K\}\}\big)\big]
.\end{equation}
We want to find an upper bound of $\Upsilon$ independent of the starting point $S_{i_{N-1}+(\delta+\gep)/a^{2}}$.
The r.h.s. of \eqref{HP} is even with respect to the starting point, therefore we can consider that
$V$ is a reflected random walk. That is why it suffices to bound from above the quantities
$W(x,a)=E_{x}(\exp(A a f \sharp\ \{i \in\{1,.,s\} :
|V_{i}|\in \{0,\dots,K\}\}))$ with $x\in \mathbb{N}$. Moreover, the Markov property
implies that $W(x,a)=W(K,a)$ for every $x\geq K$, and
$W(x,a)<W(K,a)$ if $x<K$ because the random walk starting in $K$ touches necessarily in $x$ before reaching $0$.
Therefore, we can write an upper bound of $C$, i.e.,
\begin{equation*}
\textstyle{C \leq E\Big[\exp\big(A a\sum_{j\in I_{1}}f(j)\ B_{1,N}^{j}\,
-B\sqrt{\delta}\,\Sigma_{1}\,l_{1}^{N}\big)\big|F_{N-1}
\Big]\ W(K,a)}
,\end{equation*}
and since the excursion of a random walk are independent we can assert that $B_{1,N}^{j}$
is independent of $F_{N-1}$. Hence,
\begin{multline*}
\textstyle{E\Big[\exp\big(A a\sum_{j\in I_{1}}f(j)\ B_{1,N}^{j}\,
-B\sqrt{\delta}\,\Sigma_{1}\,l_{1}^{N}\big)\big|F_{N-1}
\Big]=}\\
E\Big[\exp\big(A a \textstyle{\sum_{j\in I_{1}} }f(j)\
\sharp\big\{
i \in \big\{0,\textstyle{\frac{\delta+\epsilon}{a^{2}}}\big\} : S_{i}=j\big\}
-B\sqrt{\delta}\,\Sigma_{1}\, l_{1}^{N}\big)\Big],
\end{multline*}
and \eqref{mart}
becomes $E(M_{N}|\textsl{F}_{N-1})\leq M_{N-1}\  W(K,a)/\mu$. But $W(K,a)$ tends to $1$ as
$a\downarrow 0$ and becomes smaller than $\mu$ for $a$ small enough. That is why for $a$ small enough
$(M_{N})_{N\geq 0}$ is a super-martingale.
Since the stopping time $m_{t/a^{2}}$ is bounded from above by $t/a^{2}$, we can apply a
stopping time theorem and say that $E(M_{m})\leq E(M_{1})\leq 1$. Then, to complete
the proof of
\eqref{in1}, it suffices to show that, for
$\delta, \epsilon,a$ small enough the quantity $V_{\delta,\epsilon,a}$, defined in \eqref{Vl},
is smaller than $1$.
\begin{equation}\label{Vl}
\textstyle{V_{\delta,\epsilon,a}=\mu E\Big[\exp\big(A a \sum_{j\in I_{1}} f(j)\
\sharp\big\{
i \in \big\{0,\frac{\delta+\epsilon}{a^{2}}\big\} : S_{i}=j\big\}-B\sqrt{\delta}\,\Sigma_{1}\, l_{1}^{1}\big)\Big]}
.\end{equation}
We recall that the random walk and the local time $l_{1}^{1}$ are independent. Therefore,
\begin{equation*}
\textstyle{V_{\delta,\epsilon,a}=\mu E\Big[\exp\big(A a \sum_{j\in I_{1}} f(j)\,
\sharp\big\{
i \in \big\{0,\frac{\delta+\epsilon}{a^{2}}\big\} : S_{i}=j\big\}\big)\Big]
E\Big[\exp\big(
-B\sqrt{\delta}\Sigma_{1}\, l_{1}^{1}\big)\Big]}
.\end{equation*}

By Lemma \ref{lemsi}, we know that
\begin{equation*}
\lim_{a\to 0}V_{\delta,\epsilon,a}=\mu E\big[\exp(A\sqrt{\delta+\epsilon}\Sigma_{1} \, l_{1}^{1})\big]
E\big[\exp(-B\sqrt{\delta}\Sigma_{1}\, l_{1}^{1})\big].
\end{equation*}
Since $\Sigma_{1}$ is fixed, it enters in the constants $A$ and $B$ without changing the fact that $B>A$.
For every $x$ in $\mathbb{R}$ we denote $f(x)=E(\exp(x l_{1}^{1}))$.
The law of $l_{1}^{1}$ is known (see \cite{RevYor}), and the derivative of $f$ in $0$ satisfies
$f^{'}(0)=E(l_{1}^{1})>0$. Therefore, a first order development of $f$ gives
$f(A\sqrt{\delta+\epsilon})=1+f^{'}(0)A\sqrt{\delta+\epsilon}+
o(\sqrt{\delta+\epsilon})$
and $f(-B\sqrt{\delta})=1-f^{'}(0) B \sqrt{\delta}+o(\sqrt{\delta})$. If we take
$\epsilon\leq \delta^{2}$, we obtain
\begin{equation}\label{f}
f(A\sqrt{\delta+\epsilon})f(-B\sqrt{\delta})
\leq 1+ f^{'}(0)\sqrt{\delta} (A \sqrt{1+\delta}-B)+ o(\sqrt{\delta}).
\end{equation}
Since $B>A$, the right hand side of \eqref{f} is strictly smaller than 1 for $\delta$ small enough.
For such
a $\delta$, for $\epsilon \leq \delta^{2}$ and for $\mu>1$ but small enough we obtain
$\lim_{a \to 0} V_{\delta,\epsilon,a}<1$. As a consequence,
for $a$ small enough, $V_{\delta,\epsilon,a}<1$. This completes the proof of \eqref{in1}, and therefore,
the proof of \eqref{Kai} for $x=1$.

The proof of \eqref{Kai} for $x=2$, is easier than the former one. Indeed, $\mathbb{E}\big(
\gamma_{i}^{j}\big)<0$ for every $j\in I_{2}$, and therefore, if we choose $\beta^{''}$ such that
$\beta_{2}>\beta^{''}>\beta_{4}$,  the first order development of \eqref{Fd} gives, for $a$ small enough,
$$\log\mathbb{E}\big[\exp (2 a \beta_{2}(1+\rho^{-1}) \gamma_{i}^{j})\big]
\leq 2 a \beta''(1+\rho^{-1}) \mathbb{E}(\gamma_{1}^{j}).$$
By following the scheme of the former proof (for x=1), we notice that it suffices to replace
$\{u\in\{i_{k-1}^{v}+1,i_{k}^{v}\}
: S_{u}=j\}$ by $\{u \in \{i_{k-1}^{v}+1,i_{k-1}^{v}+(\delta+\gep)/a^{2}\}: S_{u}=j\}$ in the definition
of $M_{N}$. Moreover, there is no need to introduce $\mu>1$ in the definition of $M_{N}$, which is
in this case a positive martingale. The rest of the proof is similar to the case $x=1$.

The proof of $F_{4}<<F_{3}$ is almost the same, we just exchange the role of $\beta_{1}, \beta_{2}$ and
$\beta_{3},\beta_{4}$
in the definition of $H^{II}$. Consequently, the role of $A$ and $-B$ in
\eqref{in1} are also exchanged, and, as in the former proof, Lemma \ref{lemsi} implies the result.

\subsection{Step 4}

We notice that the quantities $m,\sigma_{1},\sigma_{2},\dots
\sigma_{m},s_{1},s_{2},\dots,s_{m}$ can also be defined for a Brownian motion on the interval $[0,t]$. In fact, we
denote $\sigma_{0}=0, z_{0}=0,$ and recursively $z_{k+1}= \inf\{s>\sigma_{k}\epsilon +\delta : B_{s}=0\}$ while
$\sigma_{k+1}$ is the unique integer satisfying $z_{k+1}\in \big((\sigma_{k+1}-1)\epsilon, \sigma_{k+1}\epsilon\big]$
and $s_{k+1}=1$ if the excursion ending in $z_{k+1}$ is in the lower half-plan, $s_{k+1}=0$ otherwise.
Finally, we let $m_{t}=\inf\{k\geq 1 : z_{k}>t\}$ and $z_{m}=t$. At this stage, we want to transform
the random walk that gives
the possible trajectories of the polymer into a Brownian
motion. For that (as in \cite{BDH}), we denote by $Q$ the measure of $(m_{t/a^{2}},\sigma_{1},\sigma_{2},\dots
\sigma_{m},s_{1},s_{2},\dots,s_{m})$ associated with the random walk on $[0,t/a^{2}]$ and by $\tilde{Q}$ the measure
of $(m_{t},\sigma_{1},\sigma_{2},\dots,\sigma_{m},s_{1},s_{2},\dots,s_{m})$ associated with the Brownian motion
on $[0,t]$.

As proven in \cite{BDH} (page $1362$) $Q$ and $\tilde{Q}$ are absolutely continuous  and their Radon-Nikod\'ym
derivative satisfies
that there exists a constant $K'_{a,\epsilon,\delta}>0$ such that for every $\delta>0$
\begin{equation}\label{Q}
\lim_{\epsilon\to 0}\limsup_{a \to 0} K'_{a,\epsilon,\delta}=0\ \  \text{and} \ \
(1-K')^{m}\leq \frac{d\tilde{Q}}{dQ}
\leq (1+K')^{m}.
\end{equation}
We recall that $\chi$ is the law of the local times $(l_{1}^{1},l_{1}^{2},\dots,l_{1}^{m})$, which are independent of the
random walk
and consequently of $Q$. Moreover, $|\overline{I}_{k}|=(\sigma_{k}-\sigma_{k-1})\epsilon/a^{2}$. Hence,
the equation \eqref{H3} gives
that $a\cdot H_{t,\gep,\gd}^{(4)}(a,h,\gb)$ depends only on $(m_{t/a^{2}},\sigma_{1},\sigma_{2},\dots
\sigma_{m},s_{1},s_{2},\dots,s_{m})$ and $(l_{1}^{1},l_{1}^{2},\dots,l_{1}^{m})$. That is why, we can write
\begin{equation*}
\textstyle{F_{t,\epsilon,\delta}^{4}(a,h,\beta_{1},\beta_{2})=\mathbb{E}\Big[\frac{1}{t} \log E_{\chi \otimes Q}\big[\exp
(a H_{t,\epsilon,\delta}^{(4)}(a,h,\beta))\big]\Big]}.
\end{equation*}
At this stage, we define $F_{5}$ by replacing the random walk by a Brownian motion, namely by integrating
over $\chi \otimes \tilde{Q}$ instead of $\chi \otimes Q$. We define
\begin{equation*}
\textstyle{H_{t,\gep,\gd}^{(5)}(a,h,\gb_{1},\gb_{2})
=H_{t,\gep,\gd}^{(4)}(a,h,\gb_{1},\gb_{2})+\frac{1}{a}\log(d\tilde{Q}/dQ)}
,\end{equation*}
and therefore,
\begin{align*}
\textstyle{F}_{t,\epsilon,\delta}^{5}(a,h,\beta_{1},\beta_{2})
=\mathbb{E}\Big[\frac{1}{t} \log E_{\chi \otimes \tilde{Q}}\big[e^{
aH_{t,\gep,\gd}^{(4)}(a,h,\gb_{1},\beta_{2})}\big]\Big]
=\mathbb{E}\Big[\frac{1}{t} \log E_{\chi\otimes Q}\big[e^{
aH_{t,\gep,\gd}^{(5)}(a,h,\gb_{1},\beta_{2})}\big]\Big].
\end{align*}
Now, we aim at proving that $F^{4}<<F^{5}$. To that aim, we calculate $H^{II}$, i.e.,
\begin{align*}
H^{II}&=H_{t,\gep,\gd}^{(4)}(a,h,\gb_{1},\beta_{2})-H_{t(1+\rho)^{2},\gep(1+\rho)^{2},\gd(1+\rho)^{2}}^{(5)}
(a(1+\rho),h^{'},\gb_{3},\beta_{4})\\
&=-\textstyle{\frac{2}{a}}(h-(1+\rho)h^{'})\sum_{k=1}^{m}s_{k}(\sigma_{k}-\sigma_{k-1})\epsilon\\
&\hspace{2.5cm}+((\beta_{1}-\beta_{3})
\Sigma_{1}+(\beta_{2}-\beta_{4})\Sigma_{2})
\textstyle{\frac{\sqrt{\delta}}{a}} \sum_{k=1}^{m} l_{1}^{k}
-\frac{1}{a(1+\rho)}\log\frac{d\tilde{Q}}{dQ}\\
&\leq -\textstyle{\frac{2}{a}}(h-(1+\rho)h^{'})\sum_{k=1}^{m}s_{k}
(\sigma_{k}-\sigma_{k-1})\epsilon-\frac{1}{a(1+\rho)}
\log\frac{d\tilde{Q}}{dQ}.
\end{align*}
We do not give the details of the end of this step because it is done in \cite{BDH} (page $1361-1362$).
To prove that $F_{5}<<F_{4}$, we consider the density
$dQ/d\tilde{Q}$ in $H^{II}$, and \eqref{Q} can also be applied. It completes the Step 4.

\subsection{Step 5}

From now on, we integrate over $\chi \otimes \tilde{Q}$ in $F^{5}$ and consequently the term
$\log\big(d\tilde{Q}/dQ\big)$ does not appear
in $H^{(5)}$ any more. In this step, transform the local times $(l_{1}^{1},\dots,l_{1}^{k},\dots)$
into the local times of the Brownian motion that determines $\tilde{Q}$.
We recall that $L_{t}$ is the local time spent at $0$ by $(B_{s})_{s\geq 0}$ between the times $0$ and $t$.

But before, we define $(R_{s})_{s\geq 0}$ a Brownian motion, independent of $B$, and
we emphasize the fact that, for every $k\in\{1,\dots,m\}$,
\begin{equation}
a \sum_{i\in \overline{I}_{k}} \hat{w}_{i}\stackrel{D}{=}R_{\sigma_{k} \gep}-R_{\sigma_{k-1} \gep}
\ \ \ \ \text{and}\ \ \ a^{2} |\overline{I}_{k}|=(\sigma_{k}-\sigma_{k-1})\gep
.\end{equation}
Then, we can rewrite the fifth Hamiltonian as
%\begin{multline}
%\textstyle{H_{t,\gep,\gd}^{(5)}(a,h,\gb_{1},\gb_{2})=-\frac{2}{a}\sum_{k=1}^{m_{t}}s_{k}[R_{\sigma_{k}\gep}
%-R_{\sigma_{k-1}\gep}
%+ h (\sigma_{k}-\sigma_{k-1})\epsilon]+
%\frac{\gb_{1}\Sigma_{1}+\gb_{2}\Sigma_{2}}{a}}\sqrt{\delta}\sum_{k=1}^{m} l_{1}^{k}
%.\end{multline}
\begin{equation}
\textstyle{H_{t,\gep,\gd}^{(5)}(a,h,\gb_{1},\gb_{2})=-\frac{2}{a}\sum_{k=1}^{m_{t}}\Big[s_{k}(R_{\sigma_{k}\gep}
-R_{\sigma_{k-1}\gep}
+ h (\sigma_{k}-\sigma_{k-1})\epsilon)-
\frac{\gb_{1}\Sigma_{1}+\gb_{2}\Sigma_{2}}{2}\sqrt{\delta}l_{1}^{k}\Big]}
.\end{equation}
We define the sixth Hamiltonian
as,
\begin{equation*}
H_{t,\gep,\gd}^{(6)}(a,h,\gb_{1},\gb_{2})=-\textstyle{\frac{2}{a}}\sum_{k=1}^{m_{t}}\Big[s_{k}
\left(R_{\sigma_{k}\gep}-R_{\sigma_{k-1}\gep}
+ h (\sigma_{k}-\sigma_{k-1}) \epsilon\right)\Big]
+\frac{\beta_{1}\Sigma_{1}+\gb_{2}\Sigma_{2}}{a}L_{t}.
\end{equation*}
At this stage, we notice that $F^{5}$ and $F^{6}$
do not depend on $a$ anymore. Hence, to simplify the following steps, we transform a bit
the general scheme of the proof. In fact, from now on, we will denote, for $i=5,6$ or $7$,
\begin{equation}
\textstyle{F_{t,\gep,\delta}^{i}(h,\beta_{1},\gb_{2})=\tilde{\mathbb{E}}\Big[\frac{1}{t} \log E_{\tilde{Q}}\big
[\exp(\overline{H}_{t,\gep,\delta}^{i}(h,\beta_{1},\gb_{2}))
\big]\Big]}\ \
\end{equation}
with $\overline{H}_{t,\gep,\delta}^{i}(h,\beta_{1},\gb_{2})=
a H_{t,\gep,\delta}^{i}(h,\beta_{1},\gb_{2})$.
Therefore, to prove that $F^{i}<<F^{j}$ we use
\begin{equation}\label{eqcha}
H^{II}=\overline{H}_{t,\gep,\delta}^{i}(h,\beta_{1},\gb_{2})-
\textstyle{\frac{1}{1+\rho}}\overline{H}_{t(1+\rho)^{2},\gep(1+\rho)^{2},
\delta(1+\rho)^{2}}^{j}
(h^{'},\beta_{3},\gb_{4})
,\end{equation}
and we show that
%\begin{equation}
$\limsup_{t \to \infty} 1/t \log \tilde{\mathbb{E}}E(\exp((1+\rho^{-1})
H^{II})))\leq 0
.$%\end{equation}

We want to prove that $F^{5}<<F^{6}$ but, by the scaling property of Brownian motion,
%\Big(namely $(B_{s})_{\{s\geq 0\}}\stackrel{D}{=}\left(\left(1+\rho\right)B_{s/(1+\rho)^{2}}\right)_
%{\{s\geq 0\}}$
%\Big)
it is not difficult to show that for $i=5$ or $6$
\begin{equation}
\overline{H}_{t(1+\rho)^{2},\gep(1+\rho)^{2},\delta(1+\rho)^{2}}^{i}(h,\beta_{1},\gb_{2})=(1+\rho)
\overline{H}_{t,\gep,\delta}^{i}((1+\rho)h,\beta_{1},\gb_{2}).
\end{equation}
Therefore, by \eqref{eqcha}, we can write $H^{II}=\overline{H}_{t,\gep,\delta}^{5}(h,\beta_{1},\beta_{2})
-\overline{H}_{t,\gep,\delta}^{6}((1+\rho)h^{'},\beta_{3},\beta_{4})$.
Thus, since $(1+\rho)h^{'}<h$ and $-\sum_{k=1}^{m}s_{k}(\sigma_{k}-\sigma_{k-1})\epsilon<0$,
we obtain
\begin{align*}
%H^{II}=-2 \ &(h-(1+\rho)h^{'}) \sum_{k=1}^{m}s_{k}(\sigma_{k}-\sigma_{k-1})\epsilon
%\ +\ \beta_{1}\Sigma_{1} \sqrt{\delta}\ \sum_{k=1}^{m}l_{1}^{k}\\
%&-\beta_{3} \Sigma_{1} \sum_{k=1}^{m} L_{p_{k}}-L_{p_{k-1}}
%+\beta_{2}\Sigma_{2} \sqrt{\delta}\ \sum_{k=1}^{m}l_{1}^{k}-
%\beta_{4} \Sigma_{2} \sum_{k=1}^{m} L_{p_{k}}-L_{p_{k-1}}\\
H^{II}\leq  \beta_{1} & \Sigma_{1}\sqrt{\delta}\ \sum_{k=1}^{m} l_{1}^{k}-
\beta_{3} \Sigma_{1}\ \sum_{k=1}^{m} L_{z_{k}^{v}}-L_{z_{k-1}^{v}}\\
&+\beta_{3} \Sigma_{1}\
(L_{t+\delta}-L_{t})
+\beta_{2}\Sigma_{2} \sqrt{\delta}\ \sum_{k=1}^{m}l_{1}^{k}-
\beta_{4} \Sigma_{2} \sum_{k=1}^{m} L_{z_{k}}-L_{z_{k-1}}
\end{align*}
with $z_{j}^{v}=z_{j}$ for every $j<m$ and $z_{m}^{v}=\inf\{t>\sigma_{m-1} \epsilon+ \delta : B_{t}=0\}$. Finally,
by the H\"older inequality,
it suffices to prove, for $B>A$, that
\begin{equation}\label{eq3}
\textstyle{\limsup_{t\to \infty} \frac{1}{t}\log E\Big[\exp\big(A \sum_{k=1}^{m} \sqrt{\delta} l^{k}_{1}
-B \sum_{k=1}^{m} L_{z_{k}^{v}}-L_{z_{k-1}^{v}}\big)\Big]\leq 0}
,\end{equation}
and
\begin{equation}\label{eq50}
\textstyle{\limsup_{t\to \infty} \frac{1}{t}\log E\Big[\exp\big(A \sum_{k=1}^{m} L_{z_{k}^{v}}-L_{z_{k-1}^{v}}
-B \sum_{k=1}^{m} \sqrt{\delta} l^{k}_{1}\big)\Big]\leq 0}
,\end{equation}
and
\begin{equation}\label{eq4}
\textstyle{\limsup_{t\to \infty}
\frac{1}{t}\log E[\exp(B(L_{t+\delta}-L_{t}))]=0}.
\end{equation}
We denote by $C_{t}$ the first time of return to the origin after time $t$.
Proving \eqref{eq4} is immediate because $C_{t}$ is a stopping time with respect to the natural filtration of $B$,
we can therefore apply the strong Markov property to obtain, for every $u\in [t,t+\delta]$, the equality
$E(\exp(B(L_{t+\delta}-L_{u}))\big|C_{t}=u)=
E[\exp(B L_{t+\delta-u})]$.
Thus, we can write
%\begin{equation*}
%E\big(\exp\left(B\left(L_{t+\delta}-L_{t}\right)\right)\big)=
%\int_{t}^{t+\delta}E\left(\exp\left(B\left(L_{t+\delta}-L_{u}\right)\right)\big|C_{t}=u\right)dC_{t}(u).
%\end{equation*}
%since  we can apply the strong Markov
%property, and we obtain
\begin{equation}\label{eq6}
E\big[\exp(B(L_{t+\delta}-L_{t}))\big]=
\int_{t}^{t+\delta}E\big[\exp(B L_{t+\delta-u})\big]dC_{t}(u)\leq
E\big[\exp(B L_{\delta})\big].
\end{equation}
This implies \eqref{eq4}, and it
remains to prove \eqref{eq3}, and \eqref{eq50}. We define a new filtration,
$F_{N}=\sigma\big(\sigma((B_{s})
_{s\leq z_{N}^{v}})\bigcup \sigma(l_{1}^{1},\dots,l^{N}_{1})\big)$. We notice that
$(z_{N}^{v})_{N\geq 0}$
is a sequence of increasing stopping times, and consequently, $F_{N}$ is an increasing filtration. We denote
by $M_{N}$ the quantity
\begin{equation}\label{M}
M_{N}=\frac{\exp\big(A \sum_{k=1}^{N} \sqrt{\delta} l^{k}_{1}-B \sum_{k=1}^{N}
L_{z_{k}^{v}}-L_{z_{k-1}^{v}}\big)}{E\big[\exp(-BL_{\delta}+A \sqrt{\delta} l^{1}_{1})\big]^{N}},
\end{equation}
which is a super-martingale with respect to $F_{N}$. Effectively,
$L$ and $(l^{k}_{1})_{k\geq 1}$ are independent, $(L_{s}-L_{z_{N}^{v}})_{s\geq z_{N}^{v}}$ is
independent of $F_{N}$ (because $B_{z_{N}^{v}}=0$) and
$L_{z_{N+1}^{v}}-L_{z_{N}^{v}}\geq L_{z_{N}^{v}+\delta}-L_{z_{N}^{v}}$. Thus, since
$E\big(\exp\big(-B(L_{z_{N}^{v}+\delta}-L_{z_{N}^{v}})\big)\big)=
E(\exp(-B(L_{\delta})))$, we obtain $E(M_{N+1}|F_{N})\leq M_{N}$.
Moreover, $m_{t}$ is a stoping time
with respect to $F_{N}$ and is bounded from above by $t/\delta$. Therefore, to prove \eqref{eq3},
it suffices to show (as in Step 3) that
for $B>A$ and $\delta$ small enough,
$V=E[\exp(A\sqrt{\delta}l_{1}^{1}-B L_{\delta})]\leq 1$.
Moreover, $L_{\delta}$ and $\sqrt{\delta}l^{1}_{1}$ have the same law and are independent. That is why we can write $V=
E[\exp(A\sqrt{\delta}l_{1}^{1})]E[\exp(-B\sqrt{\delta}l_{1}^{1})],$
which is strictly smaller than $1$ for $\delta$ small enough (as proven in Step 3).

We prove \eqref{eq50} in a very similar way. Effectively, since
$L_{z_{N+1}^{v}}-L_{z_{N}^{v}}\leq L_{z_{N}^{v}+\delta+\gep}-L_{z_{N}^{v}}$,
we prove that the inequality \eqref{eq3} is still satisfied when
$A$ and $-B$ are exchanged. Therefore, the proof of $F^{5}<<F^{6}$ is completed.
To end this step, we notice that \eqref{eq50} and \eqref{eq3} imply directly that $F^{6}<<F^{5}$.
Thus, the proof of Step 5 is completed.

\subsection{Step 6}
Let $\mu_{1}=\beta_{1}\Sigma_{1}+\beta_{2}\Sigma_{2}$ and $\mu_{3}=\beta_{3}\Sigma_{1}+\beta_{4}\Sigma_{2}$.
This step is the last one, therefore, the following Hamiltonian is the one of the continuous model,
i.e.,
\begin{equation*}
\overline{H}_{t,\gep,\gd}^{(7)}(h,\gb_{1},\gb_{2})=-2 \int_{0}^{t}\ind_{\{B_{s}<0\}}(d R_{s}+h ds)
+\mu_{1} L_{t}.
\end{equation*}
For simplicity, we define $(\phi_{s})_{s\in [0,t]}$ by $\phi_{s}=s_{k}$ for every
$s\in (\sigma_{k-1}\epsilon, \sigma_{k} \gep]$. In that way,
$\sum_{k=1}^{m} s_{k} (R_{\sigma_{k\gep}}-R_{\sigma_{(k-1)\gep}}+h(\sigma_{k}-\sigma_{k-1})\epsilon)=
\int_{0}^{t} \phi_{s} (dR_{s}+hds)$.
Moreover, the scaling property of Brownian motion gives, for $i=6$ or $7$, $$\overline{H}_{t(1+\rho)^{2},
\gep (1+\rho)^{2},\gd (1+\rho)^{2}}^{(i)}(h,\gb_1,\gb_2)\stackrel{D}{=}(1+\rho)\overline{H}_{t,\gep,\gd}^{(i)}
((1+\rho)h,\gb_{1},\gb_{2}).$$  Hence, to show that $F^{6}<<F^{7}$,
we consider (as in Step 5)
the difference $$H^{II}=\overline{H}_{t,\gep,\gd}^{(6)}(h,\gb_{1},\gb_{2})-\textstyle{\frac{1}{1+\rho}}\
\overline{H}_{t(1+\rho)^{2},\gep (1+\rho)^{2},\gd (1+\rho)^{2}}^{(7)}(h^{'},\gb_{3},\gb_{4}),
$$ which is equal to
$\overline{H}_{t,\gep,\gd}^{(6)}(h,\gb_{1},\gb_{2})-\overline{H}_{t,\gep,\gd}^{(7)}((1+\rho)h^{'},
\gb_{3},\gb_{4})$. Thus, we can bound
$H^{II}$ from above as follows
\begin{align*}
H^{II}&=-2 \int_{0}^{t} \left(\phi_{s}-\ind_{\{B_{s}<0\}}\right) dR_{s}-2\int_{0}^{t}
\left(h\phi_{s}-(1+\rho)h'\ind_{\{B_{s}<0\}}\right)ds+(\mu_{1}-\mu_{3}) L_{t}\\
H^{II}&\leq-2 \int_{0}^{t} \left(\phi_{s}-\ind_{\{B_{s}<0\}}\right) dR_{s}-2h\int_{0}^{t}
\left(\phi_{s}-\ind_{\{B_{s}<0\}}\right)ds+(\mu_{1}-\mu_{3})L_{t}
.\end{align*}
We want to prove that
%\begin{equation*}
$\limsup_{t \to \infty} \frac{1}{t} \log \tilde{\mathbb{E}}E(\exp((1+\rho^{-1})
H^{II}))
\leq 0$
%\end{equation*}
and after the integration over $\tilde{\mathbb{E}}$, it remains to prove that for $A>0$ and $B>0$ and
for $\delta,\gep$ small
\begin{equation}\label{aa}
\limsup_{t \to \infty} \textstyle{\frac{1}{t}} \log E\big[\exp\big(A\int_{0}^{t}
\big|\phi_{s}-\ind_{\{B_{s}<0\}}\big| ds-B L_{t}\big)\big]
\leq 0
.\end{equation}

As in Step 3 (see Fig. 4), we notice that between $z_{k-1}$ and $z_{k}$,
if we find an excursion of length larger than $\delta+\epsilon$, it is necessarily the one which ends at
$z_{k}$ and gives the value of $s_{k}$. It means that, apart eventually from the very beginning
of such an excursion
(between $z_{k-1}$ and $\sigma_{k-1}\epsilon$), $s_{k}$ and $\phi_{s}$ have the same value along the excursion.
Finally, we obtain
\begin{equation}\label{ab}
\textstyle{\int_{0}^{t} |\ind_{\{B_{s}<0\}}-
\phi_{s}| ds \leq P_{0,t,\delta,\gep}+ m \gep}
,\end{equation}
where $P_{u,v,\delta,\gep}$ is the sum between $u$ and $v$
of the excursion lengths which are smaller than $\delta+\gep$. The term $m \gep$
allows us to take into account the formerly mentioned situation between $z_{k-1}$ and $\sigma_{k-1}\gep$.

%We let $$\textstyle{\tilde{H}^{II}=A\int_{0}^{t}
%\big|\phi_{s}-\ind_{\{B_{s}<0\}}\big| ds-B L_{t}}$$ in \eqref{aa}. Thus, with \eqref{ab} we
%can write $\tilde{H}^{II}\leq A P_{0,t,\delta,\gep}+ A m \gep-BL_{t}$  with $A>0$
%and $B>0$. Therefore, to complete
%the proof we must show that for $\delta$ and $\gep$ small enough the inequality $\limsup_{t\to \infty} 1/t \log E(\exp(1/(1+\rho) H^{II}))\leq 0$
%occurs.
Thus, with \eqref{ab} and the H\"older inequality, we can show that the inequality
\eqref{aa} occurs if, for $\delta,\gep$ small, we have
\begin{equation}\label{eq5}
\limsup_{t\to \infty} \frac{1}{t} \log E[\exp(A\gep m- B L_{t})]\leq 0
\ \ \text{and}\ \
%\begin{equation}\label{eq6}
\limsup_{t\to \infty} \frac{1}{t} \log E[\exp(A P_{0,t,\delta,\epsilon}
-B L_{t})]\leq 0
.\end{equation}
We begin with the proof of the first inequality of \eqref{eq5}. To that aim, we recall that,
for every $k<m$, we have
$z_{k}>z_{k-1}+\delta$. Therefore, we can write $$\textstyle{A} \gep m-B L_{t}\leq A \gep m-B
\sum_{k=1}^{m} L_{z_{k-1}+\delta}-L_{z_{k-1}}+B (L_{t+\delta}-L_{t}).$$
From the equation \eqref{eq4}
and the H\"older inequality we deduce that the term $B (L_{t+\delta}-L_{t})$ does not change the result.
For this reason
we just have to consider the quantity $1/t \log E[\exp(\sum_{k=1}^{m}A \gep-B(L_{z_{k-1}+
\delta}
-L_{z_{k-1}}))]$ when $t\uparrow\infty$. As in \eqref{M}, we define the martingale
\begin{equation}\label{M2}
M_{N}=\textstyle{\frac{1}{
(V_{\gep,\delta})^{N}}}\exp\big(\sum_{k=1}^{N}
A \gep-B (L_{z_{k-1}+\delta}-L_{z_{k-1}})\big)\ \text{with}\ V_{\gep,\delta}=E[\exp(A\gep-B L_{\delta})].
\end{equation}
Since $m$ is a stopping time bounded from above by $t/\delta$, it is sufficient to show that
$V_{\gep,\delta}<1$ for $\delta,\epsilon$ small enough.
It is the case because $E[\exp(-B L_{\delta})]<1$ for every $B>0$.
Therefore, we take $\gep$ small enough and it
completes the proof.

It remains to prove the second part of \eqref{eq5}. Notice that $P_{0,t,\delta,\gep}=\sum_{k=1}^{m}
P_{z_{k-1},z_{k},\delta,\gep}$ and that for every $k\leq m$ $P_{z_{k-1},z_{k},\delta,\epsilon}\leq 2
(\delta+\gep)$ (still because there can not be more than one excursion larger than $\delta+\gep$ between $z_{k-1}$
and $z_k$). Therefore, we obtain the following upper bound $$\textstyle{A P_{0,t,\delta,\epsilon}
-B L_{t} \leq 2A(\delta+\gep)m-B
\sum_{k=1}^{m} L_{z_{k-1}+\delta}-L_{z_{k-1}}+B (L_{t+\delta}-L_{t})}.$$
As in \eqref{eq4} the term $B (L_{t+\delta}-L_{t})$
is removed, and it remains to consider $1/t \log E[\sum_{k=1}^{m}A (\gep+\delta)
-B(L_{z_{k-1}+\delta}
-L_{z_{k-1}})]$ when $t \uparrow\infty$. To that aim, we build again the martingale
\begin{equation}
\textstyle{M_{N}}=\frac{1}{(D_{\epsilon,\delta})^{N}}\exp\big(\sum_{k=1}^{N}A (\gep+\delta)
-B(L_{z_{k-1}+\delta}-L_{z_{k-1}})\big)
\end{equation}
with $D_{\epsilon,\delta}=E[\exp(A(\delta+ \gep)-B L_{\delta})]$.
The term $m$ is a bounded stopping time, therefore, it suffices to show, for $\delta,\gep$ small enough,
that $D_{\epsilon,\delta}<1$.
To that aim, we choose $\gep\leq \delta$, and it remains to consider the quantity
$E[\exp(2A\delta-BL_{\delta})]$. Moreover,
$L_{\delta}=_{D}\sqrt{\delta} L_{1}$, and
if we denote $f(x)=E[\exp(x L_{1})]$, we can use a first order development of $f$ in $0$. It gives
$f(-B\sqrt{\delta})=1-f^{'}(0) B \sqrt{\delta}+\xi_{1}(\delta) \sqrt{\delta}$ with $f'(0)>0$ and
$\lim_{x \to 0} \xi_{1}(x)=0$. We also know that,
$\exp(2A \delta)=1+2A\delta+\xi_{2}(\delta)\delta$
with $\lim_{x \to 0} \xi_{2}(x)=0$. Hence, for $\gep\leq \delta$ and $\delta$ small enough, we obtain
$E\left(\exp\left(2A\delta-BL_{\delta}\right)\right)=\exp(2A \delta)f(-B\sqrt{\delta})<1$.
The proof of $F_{6}<<F_{5}$ is exactly the same and the Step 6 is completed.

\section{Appendix}\label{sec4}

\subsection{proof of Proposition \ref{proph}}
\begin{proof}
The computation of $\tilde{\Phi}$ is based on the fact that $\tilde{\Phi}(\beta,h)$ is equal to the quantity
$h+\lim_{t \to \infty} 1/t \log E\left(\exp\left(-2h \Gamma^{-}(t)\hspace{-0.1cm}+
\hspace{-0.1cm}\beta L_{t}^{0}\right)\right)$, where $\Gamma^{-}(t)=\int_{0}^{t}\ind_{\{B_s<0\}}ds$. When
$\beta\leq 0$
we can conclude immediately that $\tilde{\Phi}(\beta,h)=h$. Therefore, in what follows we consider
$\beta>0$. Moreover,
the joint law of $(\Gamma^{-}(t),L_{t})$ is available in \cite{Kar} and takes the value
\begin{equation}\label{density}
\textstyle{dP_{\left(\Gamma^{-}(t),L_{t}^{0}\right)}(\tau,b)=\ind_{\{0<\tau<t\}}\ind_{\{b>0\}}\frac
{b\,t \exp\left(-\frac{t\, b^{2}}{8\, \tau\, (t-\tau)}\right)}{4\, \pi\, \tau^{\frac{3}{2}}\,
(t-\tau)^{\frac{3}{2}}}\ db\  d\tau}.
\end{equation}
From now on, we will denote $R_{t}=E\left(\exp\left(-2h \Gamma^{-}(t)+\beta L_{t}^{0}\right)\right)$, and with
\eqref{density} and the new variables $s=\tau/t$ and $v=b/\sqrt{t}$, we obtain
\begin{equation}\label{rt}
\textstyle{R_{t}=\int_{0}^{\infty}\frac{v\exp\big(\beta v \sqrt{t}\big)}{4 \pi}\int_{0}^{1}
\exp(-2hst) \frac{ \exp\big(
-\frac{v^{2}}{8 s (1-s)}\big)}{s^{\frac{3}{2}}(1-s)^{\frac{3}{2}}} ds dv}
.\end{equation}
In this computation we delete the constant terms because they do not change the
limit. We can write $\int_{0}^{1}$ of \eqref{rt} as the sum of $A_{1}(t)=\int_{0}^{1/2}$ and
$A_{2}(t)=\int_{1/2}^{1}$. Then, we introduce the new variable $u=s(1-s)$ in $A_{1}(t)$ and $A_{2}(t)$, and
we obtain
\begin{equation}
\textstyle{A_{1}(t)=\int_{0}^{\frac{1}{4}} \frac{\exp\big(h(\sqrt{1-4u}-1)t-\frac{v^{2}}{8u}\big)}{u^{\frac{3}{2}}
\ \sqrt{1-4u}}\,du\ \ \ \text{and}\ \ \
A_{2}(t)=\int_{0}^{\frac{1}{4}} \frac{\exp\big(-h(\sqrt{1-4u}+1)
t-\frac{v^{2}}{8u}\big)}{u^{\frac{3}{2}}
\ \sqrt{1-4u}}\,du}.
\end{equation}
It gives immediately the inequalities $A_{1}(t)\leq A_{1}(t)+A_{2}(t) \leq 2A_{1}(t)$. Therefore,
instead of studying the convergence of $1/t\,  \log R(t)$, it suffices to consider $1/t \log S(t)$ with
$S(t)=\int_{0}^{\infty} v
\exp(\beta v \sqrt{t}) A_{1}(t)dv$. We apply the Fubini Tonnelli theorem which gives
\begin{equation}\label{St}
\textstyle S(t)=\int_{0}^{\frac{1}{4}}\frac{\exp(h t\sqrt{1-4u})}{u^{\frac{3}{2}}\ \sqrt{1-4u}}\int_{0}^{\infty}
v \exp\Big(\beta v \sqrt{t}-\frac{v^{2}}{8u}\Big)dv du\ \exp(-ht)
.\end{equation}
Thus, for every $u\in [0,1/4]$, we change the variables of the second integral of \eqref{St}. To that aim, we denote
$r=v^{2}/u$. After that, we transform the variable $u$ into $x=4u$, and we obtain
\begin{equation}\label{F}
\textstyle S(t)=\frac{1}{4}\int_{0}^{1}\frac{\exp(h t\sqrt{1-x})}{\sqrt{1-x}}\ \frac{\int_{0}^{\infty}
\exp\big(\frac{\beta \sqrt{rxt}}{2}-\frac{r}{8}\big)dr}{\sqrt{x}} dx\ \exp(-ht)
.\end{equation}
The constant factor $1/4$ can be deleted and
Thus, by considering \eqref{F},
for every $\gep>0$,
we can write the following lower bound,
\begin{align*}
\textstyle\liminf_{t\to \infty}\frac{1}{t} \log S(t)+h \geq \textstyle\liminf_{t \to \infty}
\frac{1}{t}\big[\log \int_{0}^{\gep}
\frac{\exp(ht \sqrt{1-u})}
{\sqrt{1-u}\ \sqrt{u}}du+ \log \int_{0}^{\infty} e^{-\frac{r}{8}}dr\big]\geq h\sqrt{1-\gep}.
\end{align*}
Thus, we let $\gep$ tend to $0$ and we obtain
\begin{equation}\label{unom}
\textstyle\liminf_{t\to \infty}\frac{1}{t} \log S(t)+h \geq h.
\end{equation}
But we can also bound $\liminf_{t\to \infty}\frac{1}{t} \log S(t)+h$ as follows.
The laplace method %described for example in \cite{Lapla}
allows us to find the asymptotic behavior of
$Y(x)=\int_{0}^{\infty}
\exp\left(\beta \sqrt{rxt}/2-r/8\right)dr$ when $x$ tends to $\infty$. Since $\beta>0$, it gives $Y(x)
\sim_{x\to \infty} c \sqrt{xt} \exp \left(\beta^{2} xt/2\right)$ with $c>0$ that depends on $\beta$ and
we obtain
%If we choose $c^{'}<c$ and $t$ large enough, we obtain
%\begin{equation*}
%\textstyle\liminf_{t\to \infty} \frac{1}{t} \log S(t)+h\geq \liminf_{t \to \infty} \frac{1}{t} \log \int_{\gep}^{1}
%\frac{\exp(ht \sqrt{1-x})}
%{\sqrt{1-x}\ \sqrt{x}} c^{'} \sqrt{xt}\exp\left(\frac{\beta^{2}x t}{2}\right)dx
%,\end{equation*}
\begin{equation}\label{G}
\textstyle\liminf_{t\to \infty} \frac{1}{t} \log S(t)+h\geq \liminf_{t \to \infty} \frac{1}{t} \log \int_{\gep}^{1}
\frac{\exp (h t\sqrt{1-x}+\frac{t\beta^{2}x}{2})}
{\sqrt{1-x}}dx
.\end{equation}
With the formerly mentioned laplace method, we can find the asymptotic behavior of the
integral of the r.h.s. of \eqref{G}.
As $t$ tends to $\infty$, it behaves as $d \exp\big(t \big(\frac{h^{2}}{2\beta^{2}}
+\frac{\beta^{2}}{2}\big)
\big)/\sqrt{t}$ with $d>0$.
Therefore, we obtain
\begin{equation}\label{P}
\textstyle\liminf_{t\to \infty} \frac{1}{t} \log S(t)+h\geq \frac{h^{2}}{2\beta^{2}}+\frac{\beta^{2}}{2}.
\end{equation}
Finally, \eqref{unom} and \eqref{P} give
\begin{equation}\label{P1}
\textstyle\liminf_{t\to \infty} \frac{1}{t} \log S(t)+h\geq \max\big\{\frac{h^{2}}{2\beta^{2}}+\frac{\beta^{2}}{2},h\big\}.
\end{equation}
Now, we want to show that the r.h.s. of \eqref{P1} is also an upper bound of the quantity
$\limsup_{t\to \infty} 1/t
\log S(t)+h$.
To that aim, we use the fact that $\limsup_{t\to \infty} 1/t \log S(t)+h$ is equal to the maximum of
$\limsup_{t\to
\infty}1/t \log \int_{0}^{\gep}$ and $\limsup_{t\to
\infty}1/t \log \int_{\gep}^{1}$. The same kind of estimates allows us to perform the computation. Hence,
we have
\begin{equation*}
\textstyle\lim_{t\to \infty} \frac{1}{t}\log S(t)+h=\max\left(\frac{h^{2}}{2\beta^{2}}+\frac{\beta^{2}}{2},h\right)
.\end{equation*}
Finally, $\tilde{\Phi}(\beta,h)=h+\lim_{t\to \infty} 1/t \log S(t)$, and therefore,
\begin{equation*}
\textstyle\tilde{\Phi}(h,\beta)=h\ \ \text{if}\ \ h>\beta^{2}\ \
\text{and}\ \
\tilde{\Phi}(h,\beta)=\frac{h^{2}}{2\beta^{2}}+ \frac{\beta^{2}}{2}\ \ \text{if}\ \ h\leq \beta^{2}
.\end{equation*}
\end{proof}

\section*{Acknowledgments}

I am grateful to Giambattista Giacomin for his precious help and suggestions.

%%%%%%%%%%%%%%%%%%%%%%%%%%%%%%%%%%%%%%%%%%%%%%%%%%%%%%%%%%%%%%%%%%%%%%%%%%%%%%
%%%%%%%%%%%%%%%%%%%%%%%%% The bibliography %%%%%%%%%%%%%%%%%%%%%%%%%%%%%%%%%%%
%%%%%%%%%%%%%%%%%%%%%%%%%%%%%%%%%%%%%%%%%%%%%%%%%%%%%%%%%%%%%%%%%%%%%%%%%%%%%%

\end{document}